\documentclass[11pt,a4paper]{article}
\pdfoutput=1
\usepackage{amssymb}
\usepackage{amsmath}
\usepackage{graphicx}
\usepackage{color}

\renewcommand{\rho}{u}
\renewcommand{\j}{v}

\newcommand{\tc}{\tilde{c}}
\newcommand{\tA}{\tilde{A}}

\newcommand{\dt}{\Delta t}

\newcommand{\cv}{{c}}
\newcommand{\e}{{\varepsilon}}
\newcommand{\tb}{\tilde{b}}

\def\RR{\mathbb R}
\newcommand{\be}{\begin{equation}}
\newcommand{\ee}{\end{equation}}
\newcommand{\bea}{\begin{eqnarray}}
\newcommand{\eea}{\end{eqnarray}}
\newcommand{\bean}{\begin{eqnarray*}}
\newcommand{\eean}{\end{eqnarray*}}
\def\ba{\begin{array}{l}\displaystyle}
\def\ea{\end{array}}
\def\epsi{{\varepsilon}}
\def\epsilon{{\varepsilon}}

\newtheorem{remark}{Remark}
\newtheorem{definition}{Definition}
\newtheorem{proposition}{Proposition}
\newtheorem{theorem}{Theorem}
\definecolor{PBlue}{rgb}{0.0,0.0,0.69}
\definecolor{PGreen}{rgb}{0.0,0.69,0.0}
\definecolor{PMagenta}{rgb}{0.69,0.0,0.69}
\definecolor{PRed}{rgb}{0.8,0.0,0.0}

\title{A unified IMEX Runge-Kutta approach for hyperbolic systems\\ with multiscale relaxation}
\author{S. Boscarino\thanks{Mathematics and Computer Science Department, University of
Catania, Italy ({\tt boscarino@dmi.unict.it}).}\and L.
Pareschi\thanks{Mathematics Department, University of Ferrara,
Italy ({\tt lorenzo.pareschi@unife.it}).} \and
G.Russo\thanks{Mathematics and Computer Science Department,
University of Catania, Italy ({\tt russo@dmi.unict.it}).}}
\begin{document}
\maketitle
\begin{abstract}
In this paper we consider the development of Implicit-Explicit (IMEX) Runge-Kutta schemes for hyperbolic systems with multiscale relaxation. In such systems the scaling depends on an additional parameter which modifies the nature of the asymptotic behavior which can be either hyperbolic or parabolic. Because of the multiple scalings, standard IMEX Runge-Kutta methods for hyperbolic systems with relaxation loose their efficiency and a different approach should be adopted to guarantee asymptotic preservation in stiff regimes. We show that the proposed approach is capable to capture the correct asymptotic limit of the system independently of the scaling used. Several numerical examples confirm our theoretical analysis. 

\end{abstract}

 \medskip


\tableofcontents

\section{Introduction}
\label{Intr}
Hyperbolic systems with relaxation often contain multiple space-time scales which may differ of several orders of magnitude due to the various physical parameters characterizing the model. This is the case, for example, of kinetic equations close to the hydrodynamic limit \cite{BGL, Ce, CIP}. In such regimes these systems can be more conveniently described in terms of fluid-dynamic equations, when they
are considered on a suitable space-time scale.

As a prototype system, that we will use to illustrate the subsequent theory, we consider the following 
\be
\left\{  
\begin{array}{l} 
\displaystyle  
\partial_t u + \partial_x v =0, \\
\\
\displaystyle    
\partial_t v + \frac{1}{\epsi^{2\alpha}} \partial_x p(u) = -\frac{1}{\epsi^{1+\alpha}} \left( v-f(u ) \right),\quad \alpha\in[0,1]
\\
\end{array}
\right. 
\label{I61}
\ee
where $p^\prime(u)>0$. System (\ref{I61}) is hyperbolic with two distinct
real characteristics speeds $\pm \sqrt{p^\prime(u)}/\epsi^{\alpha}$.

Note that the scaling introduced in (\ref{I61}) corresponds to
the study of the limiting behavior of the solution for the usual hyperbolic system with a singular perturbation source
\be
\left\{  
\begin{array}{l} 
\displaystyle  
\partial_\tau u + \partial_\xi V =0, \\
\\
\displaystyle    
\partial_\tau V + \partial_\xi p(u) = - \frac{1}{\epsi}\left( V-F(u ) \right),
\\
\end{array}
\right. 
\label{I31}
\ee
under the rescaling $t=\epsi^{\alpha}\tau$, $\xi=x$, $v(x,t)=V(\xi,\tau)/\epsi^{\alpha}$ and
$f(u)=F(u)/\epsi^{\alpha}$. 

For $\alpha = 0$, system \label{I61} reduces to the usual hyperbolic scaling (\ref{I31}). However, if $\alpha > 0$, we are looking at larger {microscopic} times.
 
In particular, for small values of $\epsi$, using the Chapman-Enskog expansion, the behavior of the solution to (\ref{I61}) is, at least formally, governed by the following nonlinear parabolic equation
\be
\left\{
\begin{array}{l}
\displaystyle
v= f(u) - \epsi^{1-\alpha} \partial_x p(u)+\epsi^{1+\alpha}f'(u)^2 \partial_x u + {\cal{O}}(\epsi^2), \\
\\
\displaystyle
\partial_t u + \partial_x f(u )= \epsi^{1+\alpha} \partial_{x}\Bigg[\left(\frac{p'(u)}{\epsi^{2\alpha}}-f'(u)^2\right)\partial_x u\Bigg]+{\cal{O}}(\epsi^2). 
\\
\end{array}
\right.
\label{I7b}
\ee

Therefore, as $\epsilon\to 0$ when $\alpha \in [0,1)$ we obtain the scalar conservation law
\be
\left\{
\begin{array}{l}
\displaystyle
v= f(u ), \\
\\
\displaystyle
\partial_t u + \partial_x f(u )= 0. 
\\
\end{array}
\right.
\label{I7}
\ee
Note that, the main stability condition \cite{CLL, Liu} for system (\ref{I7b}) corresponds to
\be
 f^\prime (u ) ^2 < \frac{p^\prime(u)}{\epsi^{2\alpha}}, 
\label{I5}
\ee
and it is always satisfied in the limit $\epsi\rightarrow 0$ when $\alpha > 0$, whereas for $\alpha=0$ it requires suitable assumptions on the functions $f(u)$ and $p(u)$. 

In classical kinetic theory the space-time scaling just discussed leads the so-called hydrodynamical limits of the Boltzmann equation (see \cite{CIP}, Chapter 11). For $\alpha=0$ this corresponds to the \emph{compressible Euler limit}, whereas for $\alpha \in (0,1)$ the \emph{incompressible Euler limit} is obtained.     

Something special happens when $\alpha=1$. In this case, in fact, to leading order in $\varepsilon$, we obtain the parabolic equation 
\be
\left\{  
\begin{array}{l} 
\displaystyle  
v= f(u ) -\partial_{x} p(u ),\\
\\
\displaystyle    
\partial_t u + \partial_x f(u )= \partial_{xx} p(u ). 
\\
\end{array}
\right. 
\label{I4}
\ee
In other words, considering
larger times than those typical for Euler dynamics ($\alpha=1$ instead
of $\alpha \in [0,1)$), dissipative effects become non-negligible. This behavior characterizes the \emph{incompressible Navier-Stokes limit} in classical kinetic theory. 

The development of numerical methods to solve hyperbolic systems with stiff source terms in the case $\alpha=0$ has been an active area of research in the past three decades \cite{CJR, GPT,  Jin, GPT, Pa, Pe, RoAr, PR, BR}. Another series of works is concerned with the construction of robust schemes for $\alpha=1$ when a diffusion limit is obtained \cite{Klar, JL, JPT, JPT2}. However, very few papers have considered the general multiscale problem of type (\ref{I61}) for the various possible values of $\alpha$ \cite{NP,JP}. 

The common goal of this general class methods, often referred to as \emph{asymptotic-preserving} (AP) schemes, was to obtain the macroscopic behavior described by the equilibrium system by solving the original relaxation system (\ref{I61}) with coarse grids $\Delta t,~\Delta x \gg \epsi$, where $\Delta t$ and $\Delta x$ are respectively the time step and the mesh size.

Note that, since the characteristic speeds of the hyperbolic part of system (\ref{I61}) are of order $1/\epsi^{\alpha}$, most of the popular methods \cite{CJR, Jin, PR}, for the solution to hyperbolic conservation laws with stiff relaxation  present several limitations when considering the whole range of $\alpha \in [0,1]$ and fail to capture the right behavior of the limit equilibrium equation unless the small relaxation rate is numerically resolved, leading to a stability condition of the form $\Delta t \sim \varepsilon^{\alpha}\Delta x$. 
Clearly, this \emph{hyperbolic stiffness} becomes very restrictive when $\alpha > 0$, and for $\alpha = 1$ in the parabolic regime $\epsi \ll \Delta x$ where for an explicit scheme a parabolic time step restriction of the type $\Delta t \sim \Delta x^2$ is expected.  
{A special class of IMEX schemes with explicit flux and implicit relaxation is able to deal with the parabolic relaxation ($\alpha= 1$), \cite{BR}. Such methods, however, converge to an explicit scheme for the limit parabolic equation thus requiring a penalization technique to remove the final parabolic stiffness.}

In the present paper, using a reformulation of the problem, we develop high-order IMEX Runge-Kutta schemes for a system like (\ref{I61}) in the stiff regime which work uniformly with respect to the scaling parameter $\alpha$. In the parabolic regime, $\alpha = 1$, our approach gives a scheme which is not only consistent with (\ref{I4}) without resolving the small $\varepsilon$ scale, but is also capable to avoid the \emph{parabolic stiffness} leading to the CFL condition $\Delta t \sim \Delta x^2$. In other limiting regimes, that is to say when $\alpha \in [0,1)$, the scheme maintains all the nice properties of the numerical schemes for hyperbolic conservation laws, such as the ability to capture shocks with high resolution. 

Here, although the final schemes we develop will work independently on $\epsi$, we will mainly concentrate on the study of the stiff regime for system (\ref{I61}) that is to say when $\epsi \ll 1$.  

Finally we emphasize that from a physical point of
view the problem we consider here is close in spirit to the description of the macroscopic incompressible Navier-Stokes equations of fluid-dynamics by the detailed kinetic equations \cite{BGL, Ce, CIP}. Although, for the sake of simplicity, we develop our theory for one-dimensional 2$\times$2, systems the results extend far beyond these models.

The rest of the paper is organized as follows. In Section 2 we recall some of the most popular IMEX Runge-Kutta approaches and emphasize their limited applicability when $\alpha$ ranges on the whole $[0,1]$ interval. Next, in Section 3 we introduce our new approach with the aim to avoid the stiffness induced by the characteristic speeds of system (\ref{I61}). First we present the simple first order scheme and then, using the IMEX Runge-Kutta formalism, we construct high order methods. In the case $\alpha=1$, these methods give rise to an explicit approximation of the limiting parabolic problem. In Section 4 we modify the previous schemes in order to avoid the limiting parabolic stiffness. In these schemes the diffusion term in the limiting equation is integrated implicitly.
Finally in Section 5 several numerical examples are presented showing the robustness of the present approach. Some final considerations are contained in the last section and an appendix is also included.

\section{Previous IMEX Runge-Kutta approaches}
In this section, to motivate the new approach, we recall briefly other ways to tackle stiff problems through an implicit-explicit partitioning of the differential system. 

We discretize time first, and then we discretize space on the time discrete scheme. The motivation for not adopting a method of line approach is that we can  choose the space discretization which is more suitable for each term. For simplicity of notation, in the sequel we assume that $\alpha=1$, and consider the hyperbolic-to-parabolic relaxation.  Similar conclusions are obtained for $\alpha \in (0, 1)$. Considerations on the case $\alpha=0$ are reported at the end of each subsection.

\subsection{Additive approach}

Let us consider the simple implicit-explicit Euler method applied to (\ref{I61}) based on taking the fluxes explicitly and the stiff source implicitly \cite{PR}
\be
\label{eq:SP1a}
\begin{aligned}
\frac{u^{n+1}-u^n}{\dt} & = - v_x^{n}, \\
\varepsilon^2\frac{v^{n+1}-v^n}{\dt} & = -p(u^n)_x - \left(v^{n+1} - f(u^{n+1})\right). 
\end{aligned}
\ee
Solving the second equation for $v^{n+1}$ one obtains 
\be
   v^{n+1} = \frac{\epsi^2}{\epsi^2+\dt} v^n - \frac{\dt}{\epsi^2+\dt}\left(p(u^n)_x-f(u^{n+1})\right). 
\label{eq:esa}
\ee
Making use of this relation (replacing $n$ by $n-1$) in the first equation we get 
\[
\begin{aligned}
\frac{u^{n+1}-u^n}{\dt} + \frac{\epsi^2}{\epsi^2+\dt}v_x^{n-1} & = \frac{\dt}{\epsi^2+\dt}\left(p(u^{n-1})_{xx}-f(u^n)_x\right).
\end{aligned}
\]
Therefore, in the limit $\varepsilon \to 0$, we a two levels scheme (in time) for problem (\ref{I4})
\be\label{uu_n}
\begin{aligned}
\frac{u^{n+1}-u^n}{\dt}  & = p(u^{n-1})_{xx}-f(u^n)_x.
\end{aligned}
\ee
Although Eq. (\ref{uu_n}) is a consistent time discretization of the limit convection-diffusion equation (\ref{I4}), the presence of the term $u^{n-1}$ degrades the accuracy of the first order scheme.

Furthermore, since as $\varepsilon \to 0$ the equilibrium state of Eq. (\ref{eq:esa}) \[v^{n+1} = f(u^{n+1}) - p(u^n)_x\] involves two time levels, in general additional conditions on the explicit and implicit schemes are necessary in order to obtain asymptotic preserving high order methods. We refer to \cite{BR} for more details. These drawbacks are also present for any value of $\alpha \in (0,1]$. Of course, since the additive approach has been originally designed to deal with the case $\alpha=0$ there are no problems in the regime.

\subsection{Partitioned approach}
A different way to apply the implicit-explicit Euler method to (\ref{I61}) is based on taking the first equation explicitly and the second implicitly \cite{BP,BPR,NP2}
\be
\label{eq:SP2a}
\begin{aligned}
\frac{u^{n+1}-u^n}{\dt} & = - v_x^{n}, \\
\varepsilon^2\frac{v^{n+1}-v^n}{\dt} & = -p(u^{n+1})_x - \left(v^{n+1} - f(u^{n+1})\right). 
\end{aligned}
\ee
Solving for $v^{n+1}$ the second equation we get 
\be
   v^{n+1} = \frac{\epsi^2}{\epsi^2+\dt} v^n - \frac{\dt}{\epsi^2+\dt}\left(p(u^{n+1})_x-f(u^{n+1})\right), 
\ee
which substituted into the first equation results in the two level scheme
\[
\begin{aligned}
\frac{u^{n+1}-u^n}{\dt} + \frac{\epsi^2}{\epsi^2+\dt}v_x^{n-1} & = \frac{\dt}{\epsi^2+\dt}\left(p(u^{n})_{xx}-f(u^n)_x\right).
\end{aligned}
\]
Of course, since the method has been developed specifically for the case $\alpha=1$, in the limit $\varepsilon \to 0$ it yields a consistent explicit scheme for problem (\ref{I4})
\be
\begin{aligned}
\frac{u^{n+1}-u^n}{\dt}  & = p(u^{n})_{xx}-f(u^n)_x.
\end{aligned}
\ee
It is easy to verify that this approach gives a consistent explicit scheme also for any value of $\alpha \in [0,1)$ in the hyperbolic limit.
However, for non vanishingly small values of $\varepsilon$,
 the discretization of the fluxes at different time levels poses several difficulties. For example, if we consider for simplicity $p'(u)=1$, introducing the diagonal variables $u\pm \varepsilon v$ system (\ref{eq:SP2a}) reads
\be
\label{eq:SP2ab}
\begin{aligned}
\frac{(u^{n+1}+\e v^{n+1})-(u^n+\e v^n)}{\dt} & = - \frac1{\e}(u^{n+1}+\e v^{n})_x-\frac1{\e}\left(v^{n+1} - f(u^{n+1})\right), \\
\frac{(u^{n+1}-\e v^{n+1})-(u^n-\e v^n)}{\dt} & = + \frac1{\e}(u^{n+1}-\e v^{n})_x+\frac1{\e}\left(v^{n+1} - f(u^{n+1})\right), 
\end{aligned}
\ee
and therefore the space derivatives of diagonal variables in (\ref{eq:SP2ab}) are defined as a combination of two time levels 
 and it is then not clear how to discretize the system in characteristic variables. Furthermore, most numerical methods based on conservative variables evaluate the flux at the same time level. 
We refer to \cite{BPR} for a discussion on these aspects and extensions to schemes that avoid the parabolic stiffness in the relaxation limit for $\alpha=1$.

\subsection{Hybrid additive-partitioned approach}
A method which combines the advantages of the previous approaches in the various regimes has been proposed in \cite{JPT, JPT2}. The idea is to take a convex combination of the previous schemes in such a way that we have an additive scheme in hyperbolic regimes and a partitioned one in parabolic ones. This is achieved by taking
\be
\label{eq:SPh}
\begin{aligned}
\frac{u^{n+1}-u^n}{\dt} & = - v_x^{n}, \\
\varepsilon^2\frac{v^{n+1}-v^n}{\dt} & = -\phi(\varepsilon)p(u^{n})_x-(1-\phi(\varepsilon))p(u^{n+1})_x - \left(v^{n+1} - f(u^{n+1})\right), 
\end{aligned}
\ee
where $0\leq \phi(\varepsilon)\leq 1$ is such that $\phi(\varepsilon)\approx 1$ for $\varepsilon=\mathcal{O}(1)$ and $\phi(\varepsilon)\approx 0$ when $\varepsilon \ll 1$. For example $\phi(\varepsilon)=\min\{\varepsilon^2,1\}$ or the smoother approximation $\phi(\varepsilon)=\tanh(\varepsilon^2)$ were considered in \cite{JL, JPT}.

In the limit $\varepsilon \to 0$ we have the asymptotic behavior of the partitioned approach, whereas for larger values of $\varepsilon$ we have the usual additive approach. The method can be naturally extended to the general multiscale case and provides  a consistent discretization also for any value of $\alpha \in [0,1)$. We refer to \cite{JP, NP} for further details. One of the main advantages of this approach is that it results in a convex approximation at different times of the space derivative $p(u)_x$, therefore one can use different space discretizations for the derivative appearing at time $n$ and the one at time $n+1$. Typically, at time $n$ one can choose a standard hyperbolic discretization that works for large values of $\epsilon$ whereas at time $n+1$ classical central difference schemes that works in the parabolic limit suffices to avoid a CFL condition of the type $\Delta t=\mathcal{O}(\varepsilon)$. In this sense the function $\phi(\varepsilon)$ can be also interpreted as an interpolation parameter between different fluxes in the evaluation of the space derivatives in (\ref{eq:SP2a}). The determination of the optimal expression of $\phi(\varepsilon)$ in terms of stability and accuracy is a delicate aspect which is beyond the scope of the present paper.
 
\section{A unified IMEX Runge-Kutta approach} \label{SecHyp}

In this section we present a different approach which overcomes some of the drawbacks of the above mentioned methods. 

\subsection{Description of the method}
Again let us consider initially, for simplicity of notation, the case $\alpha=1$.

We now consider the following discretization for system (\ref{I61})  
\be
\label{eq:SP1}
\begin{aligned}
\frac{u^{n+1}-u^n}{\dt} & = - v_x^{n+1}, \\
\varepsilon^2\frac{v^{n+1}-v^n}{\dt} & = -\left(p(u^n)_x + v^{n+1} - f(u^n)\right). 
\end{aligned}
\ee
Solving the second equation for $v^{n+1}$ one obtains 
\be
   v^{n+1} = \frac{\epsi^2}{\epsi^2+\dt}v^n - \frac{\dt}{\epsi^2+\dt}\left(p(u^n)_x-f(u^n)\right). 
\ee
Making use of this relation in the first equation we get 
\[
\begin{aligned}
\frac{u^{n+1}-u^n}{\dt} + \frac{\epsi^2}{\epsi^2+\dt}v_x^n & = \frac{\dt}{\epsi^2+\dt}\left(p(u^n)_{xx}-f(u^n)_x\right).
\end{aligned}
\]
Note that, at variance with all the previous approaches, the first equation now uses only two time levels and the space derivatives appear all at the same time level $n$. Therefore we can rewrite the scheme in the equivalent fully explicit form
\be
\label{eq:SP2n}
\begin{aligned}
\frac{u^{n+1}-u^n}{\dt} + \frac{\epsi^2}{\epsi^2+\dt}v_x^n + \frac{\dt}{\epsi^2+\dt} f(u^n)_x& = \frac{\dt}{\epsi^2+\dt} p(u^n)_{xx},
\\
\frac{v^{n+1}-v^n}{\dt} + \frac{1}{\epsi^2+\Delta t}p(u^n)_x & = - \frac{1}{\epsi^2+\Delta t}\left(v^{n} - f(u^n)\right). 
\end{aligned}
\ee

In particular, for small values of $\dt$ the scheme (\ref{eq:SP2n}) corresponds to the system
\be
\label{eq:SP2bis}
\begin{aligned}
u_t + \frac{\epsi^2}{\epsi^2+\dt}v_x  + \frac{\dt}{\epsi^2+\dt}f(u)_x & = \frac{\dt}{\epsi^2+\dt}p(u)_{xx} + \mathcal{O}(\Delta t),\\
v_t + \frac{1}{\epsi^2+\Delta t}p(u)_x & = - \frac{1}{\epsi^2+\Delta t}\left(v - f(u)\right) + \mathcal{O}(\Delta t),
\end{aligned}
\ee
where we only used 
\[
\displaystyle \frac{u^{n+1} - u^n}{\Delta t} = u_t + \mathcal{O}(\Delta t), \quad \frac{v^{n+1} - v^n}{\Delta t} = v_t + \mathcal{O}(\Delta t)
\]
leaving all other terms.
{
Note that the left part of system (\ref{eq:SP2bis}) is hyperbolic with characteristic speeds \footnote{The subscript $1$ in the next expression indicates $\alpha=1$.}
\be
\label{char}
\lambda^{1}_{\pm}(\Delta t,\varepsilon) = \frac{\xi_1}{2}\left( \cv \pm \sqrt{\cv^2 + \frac{4\epsi^2}{(\Delta t)^2}}\right), 
\ee
with $\xi_1 = {\Delta t}/({\varepsilon^2 + \Delta t}) \in (0,1)$, $\cv = f'(u)$ and for simplicity we have set $p'(u) = 1$.  

Now, if $\dt\to 0$ for a fixed $\epsi$, system (\ref{eq:SP2bis})
converges to the original system (\ref{I61}) for $\alpha=1$ and by (\ref{char}), the characteristics speeds converge to the usual ones, i.e. $$\lambda^1_{\pm}(0,\varepsilon) =\pm \frac{1}{\varepsilon}.$$ 

On the other hand, for a fixed $\dt$, the characteristic speeds $\lambda^1_+$ and $\lambda^1_-$ in (\ref{char}) are respectively decreasing and increasing functions of $\varepsilon$ and as $\epsi\to 0$  they converge to 
\be
\label{char_lim}
\lambda^{1}_{\pm}(\Delta t,{0}) =\frac{1}{2} \left(c\pm |c|\right).
\ee
Therefore, if we denote by $\Delta x$ the space discretization parameter, we obtain the expected hyperbolic CFL condition $\Delta t \le \Delta x/|\cv|$. 

Concerning the second order derivative on the right hand side of (\ref{eq:SP2bis}) this induces a stability restriction of type $\Delta t \sim (\Delta x)^2/\xi_1$. In particular, since the discrete system (\ref{eq:SP2n}) as $\epsi\to 0$ relaxes towards 
\be
\label{eq:SP4}
\begin{aligned}
\frac{u^{n+1}-u^n}{\dt} +  f(u^n)_x & = p(u^n)_{xx},
\end{aligned}
\ee
in such a limit we have the natural stability condition which links the time step to the square of the mesh space $\Delta t \sim (\Delta x)^2$.

Similarly, in the case $\alpha \in [0,1)$, for small values of $\dt$ the scheme applied to (\ref{I61}) corresponds to the system
\be
\label{eq:SP2bisa}
\begin{aligned}
u_t + \frac{\epsi^{1+\alpha}}{\epsi^{1+\alpha}+\dt}v_x  + \frac{\dt}{\epsi^{1+\alpha}+\dt}f(u)_x & = \frac{\dt\,\epsi^{1-\alpha}}{\epsi^{1+\alpha}+\dt}p(u)_{xx} + \mathcal{O}(\Delta t),\\
v_t + \frac{\epsi^{1-\alpha}}{\epsi^{1+\alpha}+\Delta t}p(u)_x & = - \frac{1}{\epsi^{1+\alpha}+\Delta t}\left(v - f(u)\right) + \mathcal{O}(\Delta t).
\end{aligned}
\ee
The left-hand side in (\ref{eq:SP2bisa}) now has characteristic speeds
\be
\label{char_gen}
\lambda^{\alpha}_{\pm}(\Delta t,{\varepsilon}) = \frac{\xi_\alpha}{2} \left( \cv \pm \sqrt{{\cv}^{2}  + \frac{4\epsi^2}{(\Delta t)^2}}\right),
\ee
with $\xi_\alpha = {\Delta t}/({\varepsilon^{1 + \alpha} + \Delta t})\in [0,1)$. As before, if we fix $\epsi$ and send $\Delta t \to 0$ we obtain the usual characteristic speeds 
\[
\lambda^{\alpha}_{\pm}(0,\epsi)=\pm \frac1{\epsi^\alpha}.
\]
The limit behavior of the discrete system now is given by 
\be
\label{eq:SP4b}
\begin{aligned}
\frac{u^{n+1}-u^n}{\dt} +  f(u^n)_x & = 0,
\end{aligned}
\ee
and, similarly to the analysis for $\alpha=1$, we now observe that the characteristic speed do not diverge as $\varepsilon\to 0$ since we have 
\[
\lambda^{\alpha}_{\pm}(\Delta t,0)= \frac12(c\pm |c|).
\]
Therefore, we obtain the expected hyperbolic CFL condition $\Delta t \le \Delta x/|\cv|$, coming from the hyperbolic part of the system. As before, the stability restriction coming from the parabolic term, is $\Delta t \sim \Delta x^2/ \xi_{\alpha}$.
%
%

{
\subsection{Extension to general IMEX Runge-Kutta schemes}
Now we extend the analysis just performed for the simple first order implicit-explicit Euler scheme to a general IMEX-RK scheme.

\subsubsection{Notations}
An IMEX-RK scheme can be represented with a double Butcher tableau
\begin{eqnarray}\label{ars}
{\rm Explicit:} 
\>
\begin{array}{c|c}
              \tc  & \tA  \\
              \hline   \vspace{-0.25cm} \\   
      & \tb^T 
\end{array} \qquad
{\rm Implicit:} 
\>
\begin{array}{c|c}
              c & A  \\
              \hline \vspace{-0.25cm} \\
                 & b^T 
\end{array} \qquad
\end{eqnarray}
where $\tA= (\tilde{a}_{ij})$ is an $s \times s$ lower triangular matrix with zero diagonal entries for an explicit scheme, and, since computational efficiency in most cases is of paramount importance, usually the $s \times s$ matrix $A = (a_{ij})$ for an implicit scheme is restricted to the particular class of diagonally implicit Runge-Kutta (DIRK) methods, i.e., ($a_{ij} = 0,$ for $j> i$). In fact, the use of a DIRK scheme is enough to ensure that the explicit part of the scheme term is always evaluated explicitly (see \cite{ARS}, \cite{CK}, \cite{BR2009}). 

The coefficients $\tilde{c}$ and $c$ are given by the usual
relation $\tilde{c}_i =\sum_{j = 1}^{i-1}\tilde{a_{ij}}$ , $c_i =\sum_{j = 1}^{i} a_{ij}$ and the vectors $\tb = (\tb_i)_{i = 1\cdots s}$ and $b = (b_i)_{i =1 \cdots s}$ provide the quadrature weights to combine internal stages of the RK method.

The order conditions can be derived as usual matching the Taylor expansion of the exact and numerical solution, we refer to~\cite{ARS, CK, PR} for more details. Let us mention that from a practical viewpoint, coupling conditions becomes rather severe if one is interested in very high order schemes (say higher then third).  
 
Here we recall the order conditions for IMEX-RK schemes
up to order $p = 2$ and $p = 3$, under the assumption that $\tilde{c} = c$.
\begin{eqnarray}\label{OrdCond}
\begin{array}{|l|llll|}
\hline&&&&\\[-.2cm]
\mathrm{first order\,\, (consistency)}& \tilde{b}^Te = 1, & b^Te = 1.&&\\
\mathrm{second \ order}& \tilde{b} {c} = 1/2, & b^Tc = 1/2. &&\\
\mathrm{third\ order}& \tilde{b} {c}^2 = 1/3, & b^Tc^2 = 1/3, &&\\
                     &b^T \tilde{A}c = 1/6, & \tilde{b}^T\tilde{A}c = 1/6, & b^T A c=1/6, & \tilde{b}^T A c = 1/6.\\[+.2cm]
\hline
\end{array}
\end{eqnarray}
where we denote by $e^T = (1, \cdots, 1) \in \mathbb{R}^s$, and  from the previous relaxation for $\tilde{c}$ and $c$, we have $\tilde{A}e = \tilde{c}$ and $Ae = c$.

It is useful to characterize the different IMEX schemes we consider in this paper according to the structure of the DIRK method. Following~\cite{Bosc2007} we have
\begin{definition}~ 
\begin{enumerate}
\item
We call an IMEX-RK method of \emph{type I} or \emph{type A} (see \cite{PR}) if the matrix $A \in \RR^{s \times s}$ is invertible, or equivalently $a_{ii}\neq 0$, $i=1,\ldots,s$.
\item We call an IMEX-RK method of \emph{type II} or \emph{type CK} (see \cite{CK}) if the matrix $A$ can
be written as
\be
A = \left(\begin{array}{ll} 0 & 0\\
                        a  &  \hat{A}\end{array}\right),               
\label{CK1}
\ee
with $a=(a_{21},\ldots,a_{s 1})^T\in \RR^{(s-1)}$ and the submatrix $\hat{A} \in \RR^{(s-1)
\ \times \ (s-1)}$ is invertible, or equivalently $a_{ii}\neq 0$, $i=2,\ldots,s$. In the special case $a= 0$, $w_1=0$ the
scheme is said to be of \emph{type ARS} (see \cite{ARS}) and the DIRK method is reducible to a method using $s-1$ stages.
\end{enumerate}
\end{definition}


The following definition will be also useful to characterize the properties of the method \cite{BPR, BR}.
\begin{definition}
We call an IMEX-RK method \emph{implicitly stiffly accurate (ISA)} if the corresponding DIRK method is \emph{stiffly accurate}, namely
\be
a_{si} = b_i,\quad i=1,\ldots,s.
\ee
If in addition the explicit methods satisfies
\be
\tilde{a}_{si} = \tilde{b}_i,\quad i=1,\ldots,s-1
\ee
the IMEX-RK method is said to be \emph{globally stiffly accurate (GSA)}.
\end{definition}
The definitions of ISA follows naturally from the fact that $s$-stage implicit Runge-Kutta methods for which $a_{s i}=b_i$ for $i = 1, \cdots, \nu$ are called \emph{stiffly accurate} (see \cite{HW} for details). A $\nu$-stage explicit Runge-Kutta method for which $\tilde{a}_{s i}=\tilde{b}_i$ for $i = 1, \cdots s-1$ is called FSAL (\emph{First Same As Last}, see \cite{HNW} for details). Note that FSAL methods have the advantage that they require only $s-1$ function evaluations per time step, because the last stage of step $n$ coincides with the first stage of step $n+1$.

Therefore, an IMEX-RK scheme is globally stiffly
accurate if the implicit scheme is stiffly accurate and the explicit scheme is FSAL.
We observe that this definition states also that the numerical solution of a GSA IMEX-RK scheme coincides exactly with the last internal stage of the scheme.
   
\subsubsection{The unified IMEX-RK setting}   
   
Let us consider system (\ref{I61}) and again for simplicity of notation, we set $\alpha = 1$. The unified IMEX-RK approach that generalizes first order scheme (\ref{eq:SP1}) corresponds to compute first the internal stages $U$ and $V$
\begin{equation}
\label{eq:stages1}
\begin{aligned}
U & = u^n e - \dt A V_x \\
 V & =  v^n e - \frac{\dt}{\epsi^2} \tA (p(U)_x - f(U))- \frac{\dt}{\epsi^2} A V,
\end{aligned}
\end{equation}
and then the numerical solution:
\begin{equation}
\label{eq:NumSol}
\begin{aligned}
u^{n+1} & = u^n  - \dt b^T V_x \\
 v^{n+1} & =  v^n  - \frac{\dt}{\epsi^2} \tb^T (p(U)_x - f(U))- \frac{\dt}{\epsi^2} b^T V,
\end{aligned}
\end{equation}
where $f(U)$ and $p(U)$ are the vectors with component $f(U)_i = f(U^i)$ and $p(U)_i = p(U^i)$ respectively.
Solving the second equation for $V$ in (\ref{eq:stages1}) with $\zeta = \epsi^2/\Delta t$, one obtains
\begin{eqnarray}\label{Vstage}
 V  = \left( \zeta I + A\right)^{-1}\left( \zeta v^n e - \tA (p(U)_x - f(U))\right).
\end{eqnarray}
Substituting this relation in the first equation we have the resulting IMEX-RK scheme  may be written as  
\begin{equation}
\label{eq:stages22}
\begin{aligned}
\frac{U - u^n e}{\Delta t} &+ \zeta A(\zeta I + A)^{-1} e v_x^n +  A(\zeta I + A)^{-1}\tA f(U)_x &=  A(\zeta I + A)^{-1}\tA p(U)_{xx}\\
\frac{V - v^n e}{\Delta t} &+ \frac{1}{\epsi^2}\tA p(U)_x  = -\frac{1}{\epsi^2}\left(A V - \tA f(U)\right)\\
\end{aligned}
\end{equation}
and 
\begin{equation}
\label{eq:snumSol2}
\begin{aligned}
\frac{u^{n+1} - u^n}{\Delta t} &+ \zeta b^T(\zeta I + A)^{-1} e v_x^n  + b^T(\zeta I + A)^{-1}\tA f(U)_x=  b^T(\zeta I + A)^{-1}\tA p(U)_{xx}\\
\frac{v^{n+1} - v^n}{\Delta t} &+ \frac{1}{\epsi^2}\tb^T p(U)_x   = -\frac{1}{\epsi^2}\left(b^TV -\tb^T f(U)\right) .\\
\end{aligned}
\end{equation}


Then setting $B = (\zeta I + A)^{-1}$  for small values of $\Delta t$, from (\ref{eq:snumSol2}), we get the system

\begin{equation}
\label{eq:stages22b}
\begin{aligned}
u_t &+ \zeta b^T B e v_x  + b^TB\tA f(U)_x=  b^TB\tA p(U)_{xx} + \mathcal{O}(\Delta t)\\
v_t &+ \frac{1}{\epsi^2}\tb^T p(U)_x   = -\frac{1}{\epsi^2}\left(b^TV -\tb^T f(U)\right)  +  \mathcal{O}(\Delta t).\\
\end{aligned}
\end{equation}


Now, we rewrite $f(U)_x = f'(U)U_x$ and $p(U)x = p'(U)U_x$, where $f'(U)$ and $p'(U)$ are  diagonal matricies with elements $f'(U)_{ii} = f'(U^{i})$ and $p'(U)_{ii} = p'(U^{i})$ respectively. Furthermore, for simplicity we assume $f'(U^i)=c$ and $p'(U^i) = 1$.

From the IMEX-RK stages (\ref{eq:stages22}) we have

\begin{equation}
\label{eq:stages2_bis}
U= u_n e - \varepsilon^2 ABe v_x^n -  \Delta t c A B\tA U_x  + \Delta t AB \tA U_{xx},\\
\end{equation}
and using the fact that $U = u^n e + \mathcal{O}(\Delta t)$ and $b^T e=1$ (consistency of the IMEX-RK scheme) we can write the hyperbolic part in (\ref{eq:stages22}) as
\begin{equation}
\label{eq:snumSol1_bisbis}
\begin{aligned}
u_t &+ \zeta b^T B e v_x  + c b^TB\tA e u_x= \mathcal{O}(\Delta t)\\
v_t &+ \frac{1}{\epsi^2} u_x   =   \mathcal{O}(\Delta t).\\
\end{aligned}
\end{equation}
By computing the eigenvalues of the hyperbolic part we obtain
\begin{equation}\label{eigen}
\Lambda^1_{\pm}(\Delta t,{\varepsilon}) =  \frac{1}{2}\left( c b^T B\tilde{A} e \pm \sqrt{\left(c b^T B \tilde{A} e\right)^2 + 4\frac{b^TBe}{\Delta t}}\right).
\end{equation}

Next, we show that the above characteristic speeds are limited. Here we need to assume that the scheme is GSA. In fact, 
\begin{eqnarray}
b^T B  e = b^T (\zeta I + A)^{-1} e =  \Bigg\{ \begin{array}{cc} 
\displaystyle 1, & \rm{for} \ \zeta \to 0\\
\displaystyle  \sim \frac{1}{\zeta}, & \rm{for} \  \zeta \to \infty
\end{array}
\end{eqnarray}
where for the first case using the ISA property we get $b^TA^{-1}e = e^T_s e = 1$ (here $e^T_s = (0,\cdots,0,1) \in \mathbb{R}^s$), whereas for the second case we note that for the consistency of the scheme $b^T e = 1$.
On the other hand,
\begin{eqnarray}
b^T B \tilde{A} e = b^T (\zeta I + A)^{-1}\tilde{A} e =  \Bigg\{ \begin{array}{cc} 
\displaystyle 1, & \rm{for} \ \zeta \to 0\\
\displaystyle  \sim \frac{1}{2 \zeta}, & \rm{for} \  \zeta \to \infty
\end{array}
\end{eqnarray}
where for the first case using the GSA property we get $b^TA^{-1}\tilde{A}e = e^T_s\tilde{A} e = \tilde{b}^T e = 1$, whereas for the second case we note that the quantity $b^T \tilde{A}e$ is a number and if we assume that the scheme is a second order accurate one, this gives $b^T \tilde{A}e = 1/2$, by (\ref{OrdCond}).

As a consequence we have
\be
\Lambda^1_{\pm}(\Delta t,0)=\frac{1}{2}\left( c \pm \sqrt{c^2 + \frac{4}{\Delta t}}\right),\qquad \Lambda^1_{\pm}(0,\epsi)=\pm \frac1{\epsi},
\ee
and therefore the CFL condition, in the limit $\varepsilon \to 0$, becomes 
\be
\Delta t  \le \left(\frac{|\cv|}{2} + \sqrt{\frac{\cv^2}{4} + \frac{1}{\Delta t}} \right)^{-1} \Delta x
= \frac{2\sqrt{\Delta t}}{|\cv|\sqrt{\Delta t} + \sqrt{\cv^2 \Delta t + 4}}\Delta x.
\ee
Now, if $\Delta t \ll 1$ we get from
\[
\sqrt{\Delta t}  \le \frac{2}{|\cv|\sqrt{\Delta t} + \sqrt{\cv^2 \Delta t + 4}}\Delta x, 
\]
the parabolic time step restriction $\Delta t \sim \Delta x^2$.\\




  
Finally, we prove that the scheme (\ref{eq:stages22})-(\ref{eq:snumSol2}) in the limit case $\varepsilon \to 0$ is  a consistent discretization of the limit equation (\ref{I4}), i.e. the scheme is AP.

From  (\ref{Vstage}) we get
$$
V = \zeta A^{-1}v^n e - A^{-1}(I - \zeta A^{-1})\tA (f(U) - p(U)_x)  + \mathcal{O}(\zeta^2),
$$ 
and substituting in the numerical solution $v^{n+1}$ we obtain
$$
\zeta v^{n+1} = \zeta(1 - b^TA^{-1}e) v^{n} + (b^TA^{-1}\tA -\tb^T)(U_x - f(U)) - \zeta b^T{A}^{-2}\tA (f(U)-p(U)_x) + \mathcal{O}(\zeta^2).
$$
Consistency as $\zeta \to 0$, implies $(b^TA^{-1}\tA -\tb^T) = 0$, which is satisfied if the scheme is GSA, because in this case  $b^T = e^T_sA$, $\tb^T = e^T_s\tA$, therefore 
$1 - b^TA^{-1}e = 0$ and  $b^TA^{-1}\tA -\tb^T = e^T_s \tA - \tb^T =  0$.
Then we have
$$
v^{n+1} =   -b^T{A}^{-2}\tA (f(U)-p(U)_x) = e^T_s A^{-1}\tA (f(U) - p(U)_x)
$$
with 
$$
V =  - A^{-1}\tA (f(U) - p(U)_x).
$$
Similarly from (\ref{eq:stages2}) we have for the $U$ internal stages
$$
\frac{U - u^n}{\Delta t} = - \zeta v_x^n e  +  \tA(p(U)_x - f(U))_x-\zeta A^{-1}\tA(p(U)_x-f(U)) + \mathcal{O}(\zeta^2)
$$
and for the numerical solution
$$
\frac{u^{n+1} - u^n}{\Delta t} = - \zeta b^T v_x^n e  +  b^TA^{-1}\tA(p(U)_x - f(U))_x - \zeta b^TA^{-2}\tA(p(U)_x - f(U))_x + \mathcal{O}(\zeta^2),
$$
therefore, this leads to
\begin{eqnarray}\label{UnumSol}
\begin{aligned}
U &= u^n e + \Delta t \tA(p(U)_x - f(U))_x + \mathcal{O}(\xi^2),\\
u^{n+1} &= u^{n} +  \Delta t b^TA^{-1}\tA(p(U)_x - f(U))_x+ \mathcal{O}(\xi^2).
\end{aligned}
\end{eqnarray} 
Assuming that the IMEX Runge-Kutta scheme is GSA, the term $b^TA^{-1}\tA = e^T_s \tA = \tb^T$ and, as $\varepsilon \to 0$, the scheme relaxes to the explicit one, i.e.,
\begin{eqnarray}\label{ExpSol}
\begin{aligned}
\frac{U - u^n e}{\Delta t} +  \tA f(U)_x &= \tA p(U)_{xx}\\
\frac{u^{n+1} - u^{n}}{\Delta t} + \tb^T f(U)_x &= \tb^T p(U)_{xx}.
\end{aligned}
\end{eqnarray} 

\subsubsection*{Generalization to the case $\alpha \in [0, 1).$}  In the case $\alpha \in [0,1)$ analogous computations show that the characteristic speeds of the hyperbolic part are given by
\begin{equation}\label{eigen_gener}
\Lambda^\alpha_{\pm}(\Delta t, \varepsilon) =  \frac{1}{2}\left( \cv b^T B\tilde{A}e \pm \sqrt{\left(\cv b^T B \tilde{A}e\right)^2 + 4\frac{\varepsilon^{1-\alpha}}{\Delta t}b^TBe}\right)
\end{equation}
where here $\zeta = \varepsilon^{1+\alpha}/\Delta t$. We have
\[
\Lambda^\alpha_{\pm}(\Delta t, 0)=\frac12\left(c\pm |c|\right),\qquad  \Lambda^\alpha_{\pm}(0, \epsi)=\pm\frac1{\epsi^\alpha}.
\]
A similar analysis performed in this case shows that 
 \begin{eqnarray*}
\begin{aligned}
U &= u^n e - \Delta t  \tA f(U)_x +\Delta t\varepsilon^{(1-\alpha)} \tA(p(U)_{xx}+\mathcal{O}(\xi^2),\\
u^{n+1} &= u^{n} - \Delta t \tb^T f(U)_x+ \Delta t \varepsilon^{(1-\alpha)} \tb^T p(U)_{xx} +  \mathcal{O}(\xi^2).
\end{aligned}
\end{eqnarray*} 
As $\epsi \to 0$ we get the explicit Runge-Kutta scheme for system (\ref{I7})
\begin{eqnarray}\label{ExpSolGC}
\begin{aligned}
\frac{U - u^n e}{\Delta t} +  \tA f(U)_x &= 0\\
\frac{u^{n+1} - u^{n}}{\Delta t} + \tb^T f(U)_x &= 0.\\
\end{aligned}
\end{eqnarray} 

All the previous results can be stated by the following
\begin{theorem}
If the IMEX-RK scheme (\ref{eq:stages1})-(\ref{eq:NumSol}) applied to (\ref{I61}) satisfies the GSA property then, as $\varepsilon \to 0$, it becomes the explicit RK scheme characterized by the pair $(\tilde{A},\tilde{w})$
applied to the limit convection-diffusion equation (\ref{I4}) for $\alpha=1$ and to the limit scalar conservation law (\ref{I7b}) for $\alpha\in [0,1)$.
 \end{theorem}
 
Some remarks are in order. 
\begin{remark}~
\label{Remark1}
\begin{itemize}
\item The advantage of formulating the unified IMEX-RK approach in the form (\ref{eq:stages22})-(\ref{eq:snumSol2}) is that we can now adopt different space discretizations for the various term appearing in the scheme. Typically, we use classical hyperbolic-type schemes, like WENO, for the space derivatives characterizing the hyperbolic part (which now has finite characteristic speeds) and centered discretization for the second order term characterizing the asymptotic parabolic behavior.  
\item
If the IMEX-RK scheme satisfies the GSA property the numerical solution is the same as the last stage, then $V^s = v^{n+1}$ and, observing that $\tb_s = 0$, from the second equation of (\ref{eq:snumSol2}),  we have
\[
\frac{v^{n+1}-v^n}{\Delta t} + \frac1{\epsi^2} \sum_{i=1}^{s-1} \tb_{i} p(U^j)_x = -\frac1{\epsi^2}\left(\sum_{i=1}^{s-1} (b_{i} V^i - \tb_if(U^i)) + b_s v^{n+1}\right).
\]
Solving for $v^{n+1}$, after some algebra, the equation can be written as
\[
\frac{v^{n+1}-v^n}{\Delta t} + \frac1{\epsi^2+b_{s}\Delta t} \sum_{i=1}^{s-1} \tb_{i} p(U^j)_x = -\frac1{\epsi^2+b_{s}\Delta t}\left(\sum_{i=1}^{s-1} (b_{i}V^i-\tb_{i} f(U^i)) +b_s v^n\right).
\]
Assuming $p'(u)=1$ and using the fact that $U=u^n e+ {\cal O}(\Delta t)$
we can now write the hyperbolic part in (\ref{eq:stages22}) as
\begin{equation}
\label{eq:snumSol1_bisbisg}
\begin{aligned}
u_t &+ \zeta b^T B e v_x  + c b^TB\tA e u_x= \mathcal{O}(\Delta t)\\
v_t &+ \frac{1}{\epsi^2+a_{ss}\Delta t} u_x   =   \mathcal{O}(\Delta t).\\
\end{aligned}
\end{equation} 
Therefore, the characteristic speeds read
\begin{equation}\label{eigeng}
\Lambda^1_{\pm}(\Delta t,{\varepsilon}) =  \frac{1}{2}\left( c b^T B\tilde{A} e \pm \sqrt{\left(c b^T B \tilde{A} e\right)^2 + 4\frac{\zeta b^TBe}{\epsi^2+a_{ss}\Delta t}}\right),
\end{equation}
and now when $\epsi\to 0$, as in the first order case discussed in Section \ref{SecHyp}, we get
\[
\Lambda^1_{\pm}(\Delta t,0)=\frac{1}{2}\left( c \pm |c|\right).
\]
\end{itemize}
\end{remark} 
 

%
}

\subsection{Removing the parabolic stiffness} \label{SecPar}
Although the final schemes developed in the previous section will work independently on $\varepsilon$ and $\alpha$, for small values of $\varepsilon$  they relax to an explicit RK scheme originating a time step restriction of the type $\Delta t \approx \Delta x^2/\epsi^{1-\alpha}$, for $\alpha\in (0,1]$. Therefore, for small $\varepsilon$, only the case $\alpha =1$ poses stability restriction. For this reason, we shall consider $\alpha = 1$ in this subsection.


The natural idea here is to treat the term $p(u)_x$ in (\ref{I61}) implicitly and to observe that in the limit case, i.e. $\varepsilon \to 0$, the IMEX-RK scheme relaxes to an IMEX-RK scheme for the limit convection-diffusion equation (\ref{I4}) where the diffusion term is now evaluated implicitly. 
Compared to similar schemes presented in \cite{BPR, BR, BLR} this new approach has the advantage that is not based on a penalization technique and therefore avoids the difficult problem of the optimal determination of the penalization parameter (see \cite{BPR}).


%


The unified IMEX-RK scheme for system (\ref{I61}), now reads
\begin{equation}
\label{eq:stages1bis}
\begin{aligned}
U & = u^n e - \dt A V_x \\
 V & =  v^n e + \frac{\dt}{\epsi^2} \tA f(U) - \frac{\dt}{\epsi^2} A (V+p(U)_x).
\end{aligned}
\end{equation}
and 
\begin{equation}
\label{eq:snumSol1bis}
\begin{aligned}
u^{n+1}  &=  u^n  - \Delta t \,b^T V_x\\
v^{n+1} &=  v^n    +  \frac{\dt}{\epsi^2}\tb^T f(U) - \frac{\dt}{\epsi^2}b^T(V + p(U)_x).\\
\end{aligned}
\end{equation}

Now, solving the second equation in (\ref{eq:stages1bis}) for $V$, one obtains
\begin{eqnarray}\label{Vstage2}
 V  = \left( \zeta I + A\right)^{-1}\left( \zeta v^n e + \tA f(U) - A p(U)_x \right),
\end{eqnarray}
where $\zeta=\epsilon^2/\Delta t$. Using this relation in the first equation of (\ref{eq:stages1bis}) we get for the internal stages
\begin{equation}
\label{eq:stages2}
\begin{aligned}
\frac{U - u^n}{\Delta t} + \zeta A(\zeta I + A)^{-1}v_x^n e  &=  A(\zeta I + A)^{-1}(A p(U)_x -\tA f(U))_x,\\
\frac{V - v^n}{\Delta t} + \frac{1}{\epsi^2}\tA f(U)  &= -\frac{1}{\epsi^2}A (V+p(U)_x),\\
\end{aligned}
\end{equation}
and similarly for the numerical solution
\begin{equation}
\label{eq:snumSol1}
\begin{aligned}
\frac{u^{n+1} - u^n}{\Delta t} + \zeta b^T(\zeta I + A)^{-1}v_x^n e  &=  b^T(\zeta I + A)^{-1}(Ap(U)_x - \tA f(U))_x,\\
\frac{v^{n+1} - v^n}{\Delta t} + \frac{1}{\epsi^2}\tb^T f(U)  &= -\frac{1}{\epsi^2}b^T(V + p(U)_x).\\
\end{aligned}
\end{equation}

Then setting $B = (\zeta I + A)^{-1}$  for small values of $\Delta t$, from (\ref{eq:snumSol1}), we get the system

\begin{equation}
\label{eq:stages22bb}
\begin{aligned}
u_t &+ \zeta b^T B e v_x  + b^TB\tA f(U)_x=  b^TB A p(U)_{xx} + \mathcal{O}(\Delta t)\\
v_t &+ \frac{1}{\epsi^2} b^T p(U)_x   = -\frac{1}{\epsi^2}\left(b^T V -\tb^T f(U)\right)  +  \mathcal{O}(\Delta t),\\
\end{aligned}
\end{equation}
which is similar to Eq. (\ref{eq:stages22b}), except that now $A$ and $b^T$ appear in front of $p(U)$ terms in place of $\tA$ and $\tb$ respectively. 
%
Therefore, we have the same characteristic speeds as in the unified IMEX-RK approach presented in the previous section and the same conclusions on the hyperbolic CFL condition holds true. 

Concerning the AP property, in the limit $\zeta\to 0$ one has from (\ref{Vstage2})
\[
V =  -  p(U)_x  + A^{-1}\tA f(U),
\]
and by the GSA property, i.e. $b^TA^{-1}\tA = e^T_s \tA = \tb^T$, we get from (\ref{eq:snumSol1})
$$
v^{n+1} = - e^T_s(p(U)_x - A^{-1}\tA f(U)). 
$$




Thus, scheme (\ref{eq:stages2})-(\ref{eq:snumSol1}), becomes an IMEX-RK scheme for the convection-diffusion equation (\ref{I4})
\begin{eqnarray}\label{IMEXSol}
\begin{aligned}
\frac{U - u^ne}{\Delta t}  + \tA f(U)_x  & =  A p(U)_{xx},\\
\frac{u^{n+1} - u^{n}}{\Delta t} +\tb^T f(U)_x &= b^T p(U)_{xx}.
\end{aligned}
\end{eqnarray} 

Clearly, scheme (\ref{IMEXSol}) is a consistent approximation of the limit equation (\ref{I4}) where now the diffusion term is evaluated implicitly, therefore the CFL condition of such scheme is uniquely determined by the hyperbolic restriction $\Delta t \sim \Delta x$. 


{
A similar analysis in the case $\alpha \in [0,1)$ for scheme (\ref{eq:stages1bis})-(\ref{eq:snumSol1bis}) under the GSA assumption produces the explicit Runge-Kutta scheme
\begin{eqnarray}\label{}
\begin{aligned}
U &= u^n - \Delta t \tA f(U)_x,\\ 
u^{n+1} &= u^{n} - \Delta t \tilde{b}^T\tA f(U)_x.\\
\end{aligned}
\end{eqnarray} 
Therefore we can summarize the results in the following
\begin{theorem}
If the IMEX-RK scheme (\ref{eq:stages1bis})-(\ref{eq:snumSol1bis}), applied to (\ref{I61})  for $\alpha=1$, satisfies the GSA property then, as $\varepsilon \to 0$, it becomes the IMEX-RK method characterized by the pairs $(\tilde{A},\tilde{b})$ and $(A, b)$ for the
limit convection-diffusion equation (\ref{I4}). Otherwise for $\alpha \in [0,1)$ the IMEX-RK scheme as $\varepsilon \to 0$ yields the explict RK method characterized by the pair $(\tilde{A},\tilde{b})$ for the limit scalar conservation law (\ref{I7}).
\end{theorem}
}
\section{Numerical applications}
In this section we present numerical results that confirm the validity of the new approach presented in section \ref{SecPar}. In all tests we used the last approach described in  section \ref{SecPar}, so that we remove the parabolic restriction in the limit of small $\varepsilon$ and $\alpha = 1$.
All the numerical examples presented here refer to the IMEX-RK schemes reported in the Appendix.  We shall use the notation NAME$(\nu,\sigma,p)$, where the triplet $(\nu, \sigma, p)$ where $\nu, \sigma$ and $p$ represent respectively  the number of explicit function evaluations, the number of implicit function evaluations and the order of accuracy.

In order to avoid spurious numerical oscillations arising near discontinuities of the solutions, we use interpolating non-oscillatory algorithms, like WENO method, \cite{Shu}. In these numerical test, as emphasized in Remark~\ref{Remark1}, we use classical hyperbolic-type schemes, like finite difference discretization with  WENO reconstruction, for the space derivatives characterized the hyperbolic part while for the second order term $p(u)_{xx}$ we used the standard $2$-th and $4$-th order finite difference technique. Note that when we consider third-order IMEX R-K scheme we use  $4$-th order finite difference to discretize the term $p(u)_{xx}$ except at the nearby boundary points where a $3$-rd order formula was implemented, and this guarantees  to achieve a global third-order.
%

In all our numerical results we take $\Delta t = \lambda_{CFL} \Delta x$ in all regimes. The precise choice of $\lambda_{CFL}$ is reported in the figure captions.


As a mathematical model for our numerical experiments we consider the Ruijgrook-Wu model of the discrete kinetic theory of rarefied gases. The model describes a two-speed gas in one space dimension and corresponds to the system \cite{GP, JPT, RW}
\be
\left\{
\begin{array}{l} 
\displaystyle  
M \partial_t f^+ + \partial_x f^+=-\frac{1}{\rm Kn} ( a f^+ - b f^- - c f^+ f^-),\\
\\
\displaystyle
M \partial_t f^- - \partial_x f^-=\frac{1}{\rm Kn} ( a f^+ - b f^- - c f^+ f^- ), 
\\
\end{array}
\right. 
\label{E3}
\ee
where $f^+$ and $f^-$ denote the particle density distribution at time $t$, position $x$ and with velocity $+1$ and $-1$ respectively. Here ${\rm Kn}$ is the Knudsen number, $M$ is the Mach number of the system and $a$,$b$ and $c$ are positive constants which characterize the microscopic interactions. The local (Maxwellian) equilibrium is defined by 
\be
f^+ = \frac{bf^-}{a-cf^-}.
\ee
The macroscopic variables for the model are the density $\rho$
and momentum $v$ defined by
\be
\rho = f^++f^-,\quad \j=(f^+ -f^-)/M.
\ee
The nondimensional multiscale problem is obtained taking $M=\epsi^\alpha$ and ${\rm Kn}=\epsi$,
%
%
%
%
the Reynolds number of the system is then defined as usual according to $Re=M/Kn=1/\epsi^{1-\alpha}$. The model, as we will see, has the nice feature to provide nontrivial limit behaviors for several values of $\alpha$ including the corresponding compressible Euler ($\alpha=0$) limit and the incompressible Euler ($\alpha \in (0,1)$) and Navier-Stokes ($\alpha=1$) limits.


\subsubsection*{Test 1. Diffusive scaling in the linear case}
For this numerical test we consider the case $\alpha=1$ with $c=0$, $a=1+A\varepsilon $ and $b=1-A\varepsilon$ in the r.h.s. of (\ref{E3}).   

Adding and subtracting  the two equations in (\ref{E3})
one obtains the following macroscopic equations for $\rho$ and
$\j$
 
\be
\left\{  
\begin{array}{l} 
\displaystyle  
\partial_t\rho + \partial_x \j =0, \\
\\
\displaystyle    
\partial_t \j + \frac{1}{\epsi^2} \partial_x \rho = -\frac{1}{\epsi^{2}}(\j - A \rho).
\\
\end{array}
\right. 
\label{E1}
\ee
In the limit $\varepsilon \to 0 $ the second equation relaxes to the local equilibrium 
\[
\j = A\rho - \frac{\partial \rho}{\partial x}
\]
and substituting in the fist equation this gives the limiting advection-diffusion equation
\[
\rho_t + A \rho_x = \rho_{xx}. 
\] 
Next we fix $A = 1$ and observe that the limiting advection-diffusion equation admits the 
exact solution 
\begin{eqnarray}\label{InCondBW}
\rho(x,t) = e^{-t}\sin(x-t), \quad \j(x,t) = e^{-t}(\sin(x-t) - \cos(x-t)),
\end{eqnarray}
on the domain $[-\pi, \pi]$ with periodic boundary conditions. We start to consider as initial conditions (\ref{InCondBW}) at $t = 0$, and choose $\varepsilon = 10^{-6}$ with final time $T= 0.1$.  The numerical results are compared with (\ref{InCondBW}) at $t = T$. Relative errors and convergence orders are reported in Tables \ref{tab: ConvRate1}-\ref{tab: ConvRate4}. In order to check the temporal order of convergence of these schemes, $\Delta x$ decreases with the time step $\Delta t$ accordingly to the CFL condition $\Delta t = 0.5\, \Delta x$. In the Tables we show the order of convergence as
\[ p = \log_2(||E_{\Delta t}||_{\infty}/ ||E_{\Delta t/2 }||_{\infty}),\]
with $E_{\Delta t}$ the relative error computed with time step $\Delta t$. We consider the error obtained with $N = 40$ time steps up to $N = 640$ time steps in the interval $[0,T]$.
We observe that the classical order of the methods is maintained in the limit case for the density $\rho$. We omit the convergence results of ARS(2,2,2) scheme (\ref{ars2}) since they are essentially the same as the ones obtained with CK(2,2,2) scheme reported in (\ref{PtypeCK}).

Typically, classical schemes, as ARS(2,2,2), CK(2,2,2) and BPR(3,4,3),  which satisfy only the GSA property, have few numbers of internal stages but maintain, in the limit $\varepsilon \to 0$, the order in time only for the $u$-component, while for the $v$-component they do not guarantee even the consistency \cite{BR}. 
As an example, we report the convergence table for the variable $\j$ in Table~\ref{tab: ConvRate4bis} for the CK(2,2,2) scheme (\ref{PtypeCK})  and we  observe that we do not obtain the correct classical order of accuracy for $\j$ in the limit case $\varepsilon \to 0$. 

The reason for the lack of consistency in the algebraic component $v$, is that this classical schemes does not satisfy the additional order conditions (\ref{AddOrdCond2}) (see Appendix) that guarantees the consistency of the scheme and the order up to $2$ for the variable $v$.  As a comparison we construct a new scheme, BPR(4,4,2), that satisfies the additional order conditions (\ref{AddOrdCond2}), and in Table \ref{tab: ConvRate4} we observe the correct convergence rate for both components, $\rho$ and $\j$. Furthermore, in Fig. \ref{CompV}, we also compare the numerical solutions (star points) for $\rho$ and $\j$ obtained with the second order ARS(2,2,2) scheme (\ref{ars2}) and  BPR(4,4,2) scheme. In the same figure, the exact solution is plotted by a continuous line.  Note again that the new scheme BPR(4,4,2) shows the correct behaviour for the $v$ variable.
 
Of course there are some advantage and disadvantages to consider additional order conditions (\ref{AddOrdCond2}). The advantage is that they guarantee the correct order in the limit case ($\varepsilon \to 0$) for both variables. On the other hand, to construct schemes that satisfy (\ref{AddOrdCond2}), requires a larger number of internal stages, due to these extra order conditions.
  \begin{figure}[tb]
\begin{center}
\includegraphics[height=2.9in]{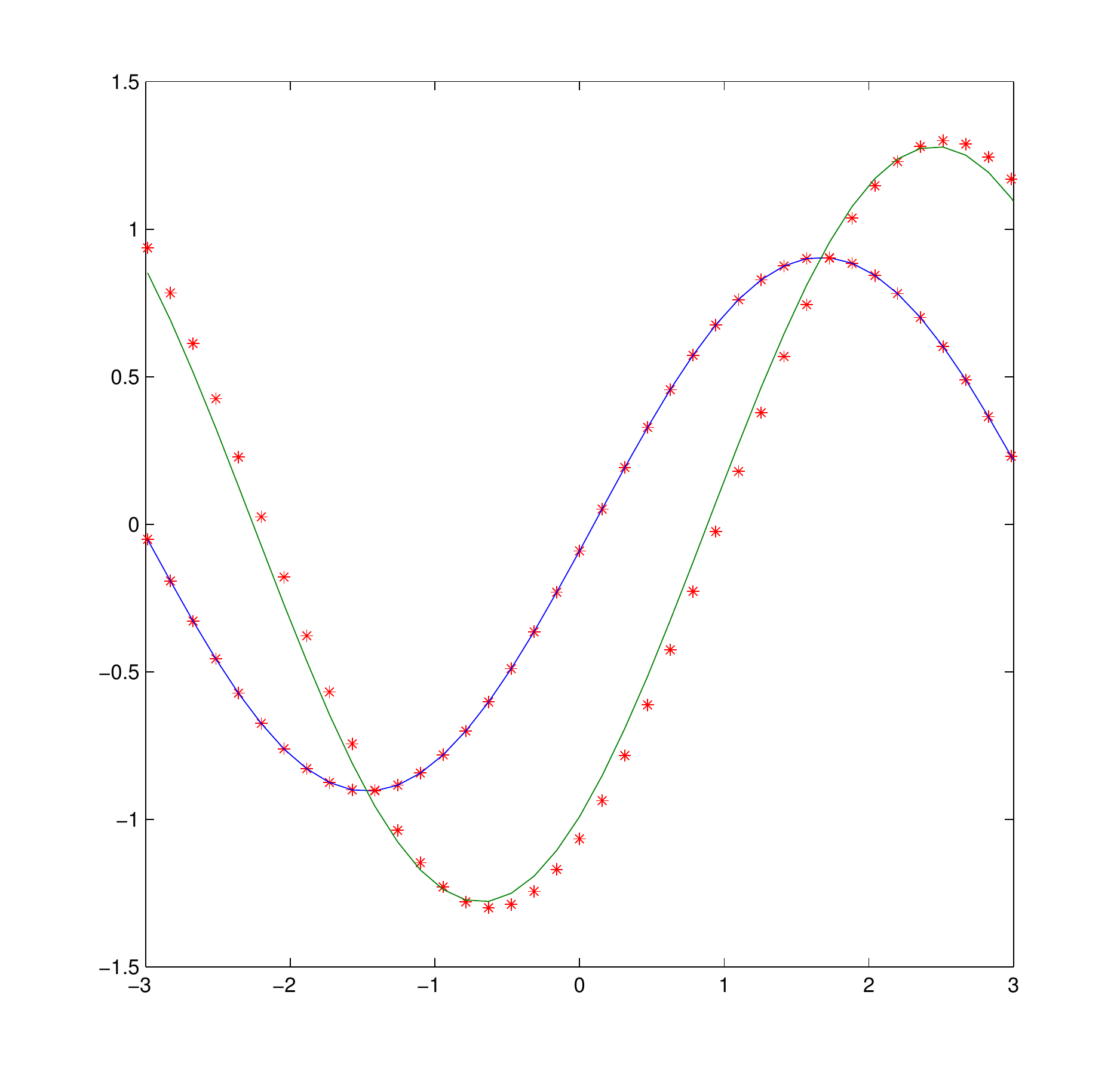},
\includegraphics[height=2.9in]{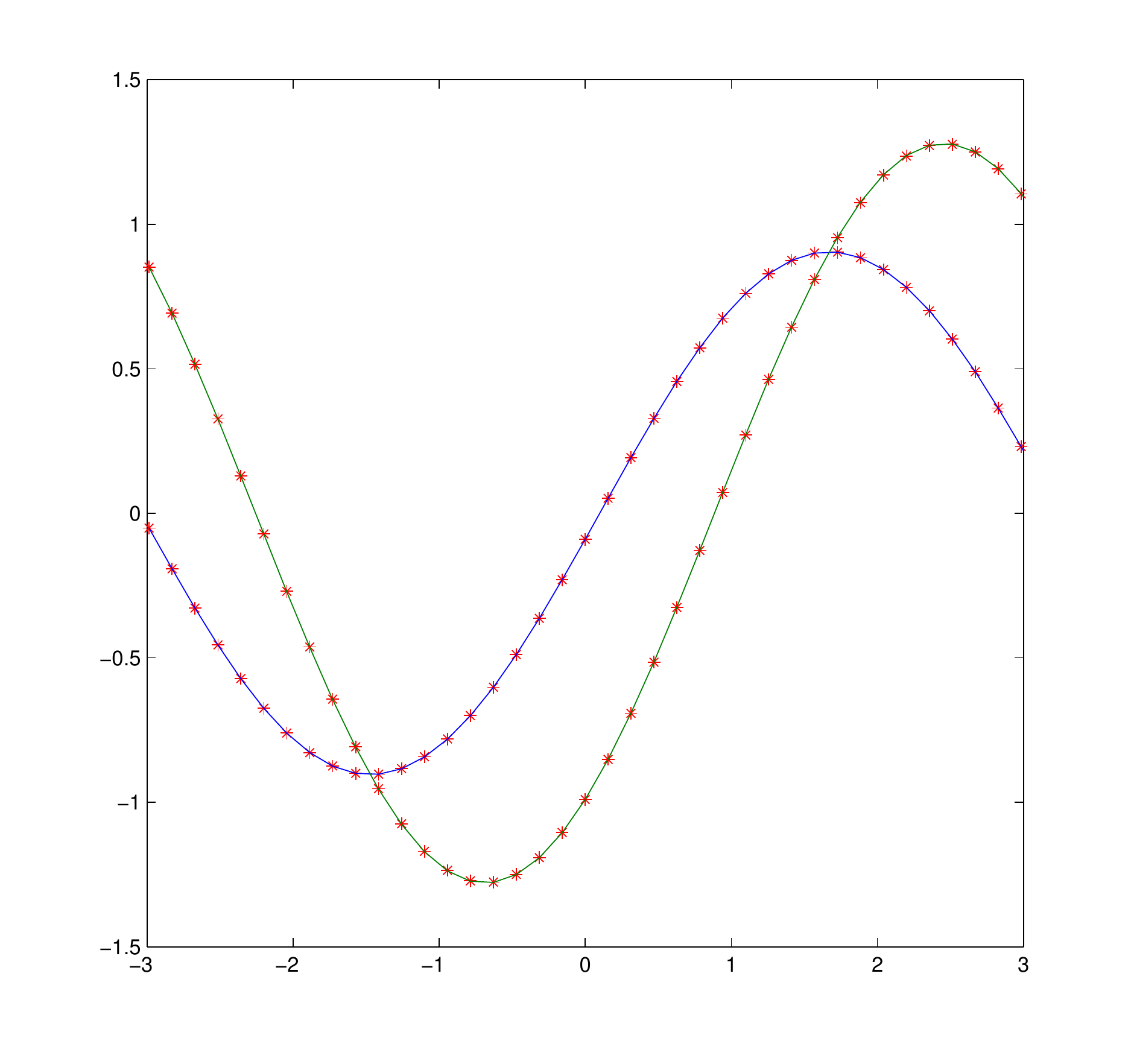}
\caption{Test 1. Comparison between classical ARS(2,2,2) scheme (left) and the BPR(4,4,2) scheme (right) with $N = 40$ and $\Delta t = \Delta x$.}
\label{CompV}
\end{center}
\end{figure}
\begin{table}
\begin{center}
\caption{{Test 1. Converge rates for the density $\rho$ with $\varepsilon = 10^{-6}$.
}\label{tab: ConvRate1}
}
\bigskip
\begin{tabular}{c | c|c|c}
\hline
\cline{1-3} Method & N& $L_{\infty}$ error &Order\\
\hline
\cline{1-3}   ARS(1,1,1) & 40&$6.4800e-03$ & $--$ \\
\hline
\cline{1-3}   ARS(1,1,1)& 80& $3.5082e-03$ & $0.8853$ \\
\hline
\cline{1-3}    ARS(1,1,1)  &160& $1.9203e-03$ & $0.8694$ \\
\hline
\cline{1-3}   ARS(1,1,1)  &320& $9.6447e-04$ & $0.9935$\\
\hline
\cline{1-3}   ARS(1,1,1)  &640& $4.8457e-04$ & $0.9930$\\
\hline
\hline
\cline{1-3}   CK(2,2,2) & 40& $1.4911e-04$ & $--$ \\
\hline
\cline{1-3}   CK(2,2,2)& 80& $3.9405e-05 $ & $1.9199$ \\
\hline
\cline{1-3}    CK(2,2,2)  & 160& $1.1356e-05 $ & $ 1.7949$ \\
\hline
\cline{1-3}   CK(2,2,2)  & 320& $2.8331e-06 $ & $ 2.0030$\\
\hline
\cline{1-3}   CK(2,2,2)  & 640& $7.0874e-07 $ & $1.9991$\\
\hline
\hline
\cline{1-3}   BPR(3,4,3) & 40& $5.8318e-06$ & $--$ \\
\hline
\cline{1-3}   BPR(3,4,3)& 80& $7.8658e-07$ & $2.8903$ \\
\hline
\cline{1-3}    BPR(3,4,3)  & 160& $1.2095e-07$ & $2.7012$ \\
\hline
\cline{1-3}   BPR(3,4,3)  & 320& $1.5297e-08$ & $2.9831$\\
\hline
\cline{1-3}   BPR(3,4,3)  & 640& $1.9253e-09$ & $2.9901$\\
\hline
\end{tabular}
\end{center}
\end{table}

\begin{table}[htb]
\begin{center}
\caption{{Test 1. Converge rates for $v$ for the scheme CK(2,2,2) with $\varepsilon = 10^{-6}$.}
}\label{tab: ConvRate4bis}
\bigskip
\begin{tabular}{c | c|c|c}
\hline
\cline{1-3} Method & N& $L_{\infty}$ $v$-error &Order of $v$  \\
\hline
\cline{1-3} CK(2,2,2) & 40 & $1.3156e-02$& $--$   \\
\hline
\cline{1-3}  CK(2,2,2) & 80 &$1.4897e-02$ & $-0.1793$\\
\hline
\cline{1-3}   CK(2,2,2) & 160& $1.1520e-03$ & $3.6928$ \\
\hline
\cline{1-3}   CK(2,2,2)  & 320& $1.2584e-03$& $-0.1274$\\
\hline
\cline{1-3}   CK(2,2,2)  &    $640$ &  $1.2877e-03$& $-0.0332$\\
\hline
\end{tabular}
\end{center}
\end{table}
\begin{table}[htb]
\begin{center}
\caption{{Test 1. Converge rates for $\rho$ and $\j$ for the scheme (\ref{New}) with $\varepsilon = 10^{-6}$.
}\label{tab: ConvRate4}
}
\bigskip
\begin{tabular}{c | c|c|c|c|c}
\hline
\cline{1-3} Method & N& $L_{\infty}$ $\rho$-error &Order $\rho$& $L_{\infty}$ $\j$-error & Order $\j$\\
\hline
\cline{1-3}   BPR(4,4,2) & 40& $ 1.9129e-04$ & $--$ & $ 2.8704e-04$ & $--$\\
\hline
\cline{1-3}   BPR(4,4,2)& 80& $ 4.9963e-05$ & $1.9368$ & $8.0261e-05$ & $ 1.8385$\\
\hline
\cline{1-3}    BPR(4,4,2)& 160& $1.4374e-05$ & $1.7974$ & $2.0603e-05$&$ 1.9618$\\
\hline
\cline{1-3}   BPR(4,4,2) & 320& $  3.5895e-06$ & $ 2.0016$ & $ 5.2702e-06$&$1.9669$\\
\hline
\cline{1-3}   BPR(4,4,2) & 640& $9.0120e-07$ & $1.9939$ &$1.4011e-06$ &$ 1.9113$\\
\hline

\end{tabular}
\end{center}
\end{table}
Next, we consider a Riemann problem with initial data
\be
\left\{
\begin{array}{ll}
\rho_L = 4.0,\quad \j_L=0,\quad & 
-10 < x < 0,
\\
\rho_R = 2.0,\quad \j_R=0,\quad &
0 < x < 10,
\\
\end{array} \right. \label{RWdata}
\ee
and inflow and outflow boundary conditions.  The exact solution for the limit advection-diffusion equation is 
\[
\rho(x,t) = \frac{1}{2}(\rho_L + \rho_R) + \frac{1}{2}(\rho_L - \rho_R){\rm erf}\left(\frac{t-x}{2\sqrt{t}}\right),
\]
where {erf}$(x)$ denotes the error function. We take $\Delta x = 0.2$ and  final time $T = 3.0$. In Figure \ref{LinearCase} we compare the numerical solutions computed by schemes SP(1,1,1), BPR(2,4,4) and BPR(3,3,5)  for the mass density in the intermediate regime ($\varepsilon = 0.5$) and in diffusive one ($\varepsilon = 10^{-6}$)  with a {reference} solution obtained with $\Delta x = 0.001$. As we can see, the numerical results for the different schemes describe the exact motion of the shock and, in the small relaxation limit, are in excellent agreement with the analytical ones using a hyperbolic time step $\Delta t = 0.5 \Delta x$. 
\begin{figure}[htbp]
\begin{center}
\centering
\includegraphics[width=3.1in]{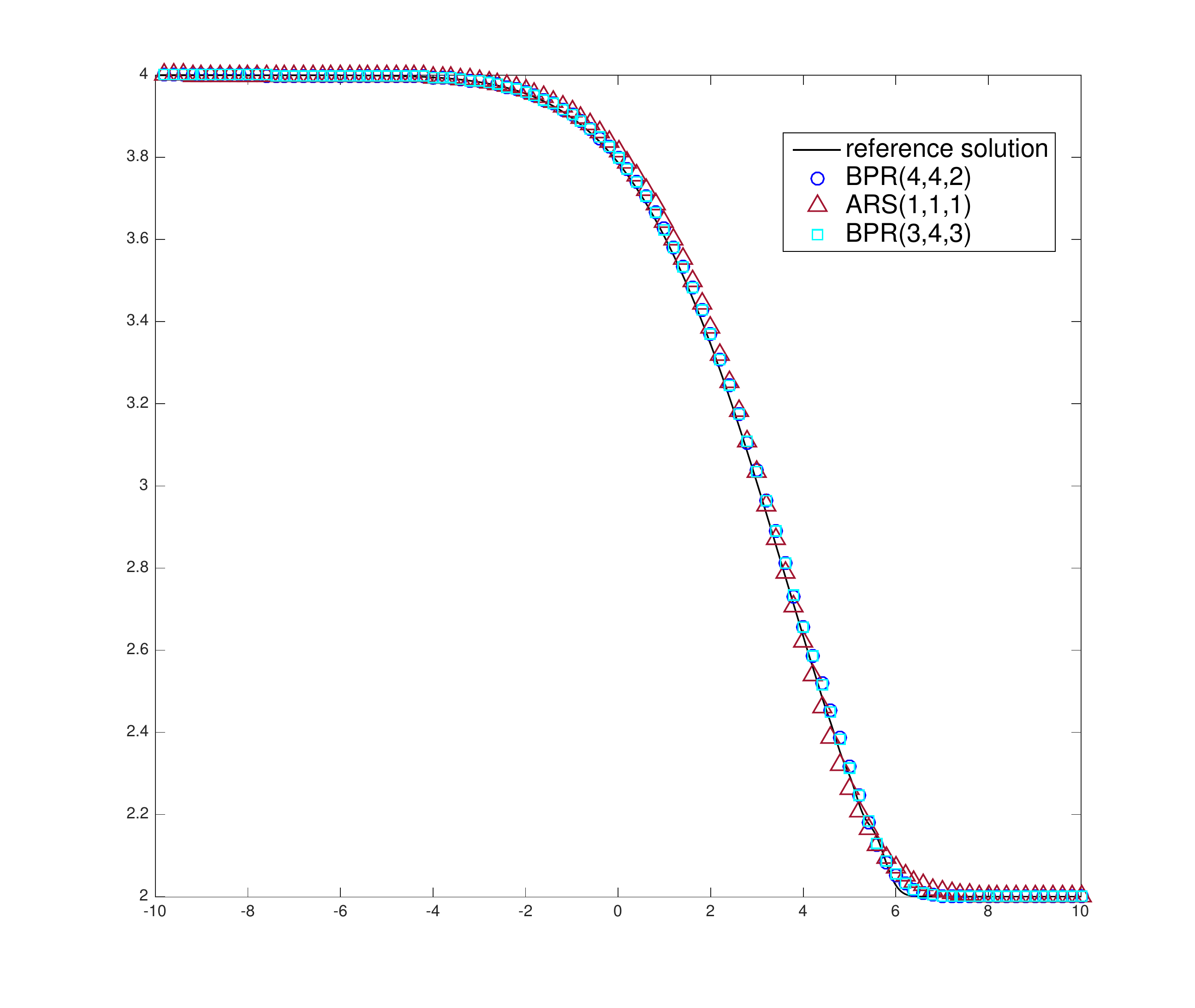}
\includegraphics[width=3.1in]{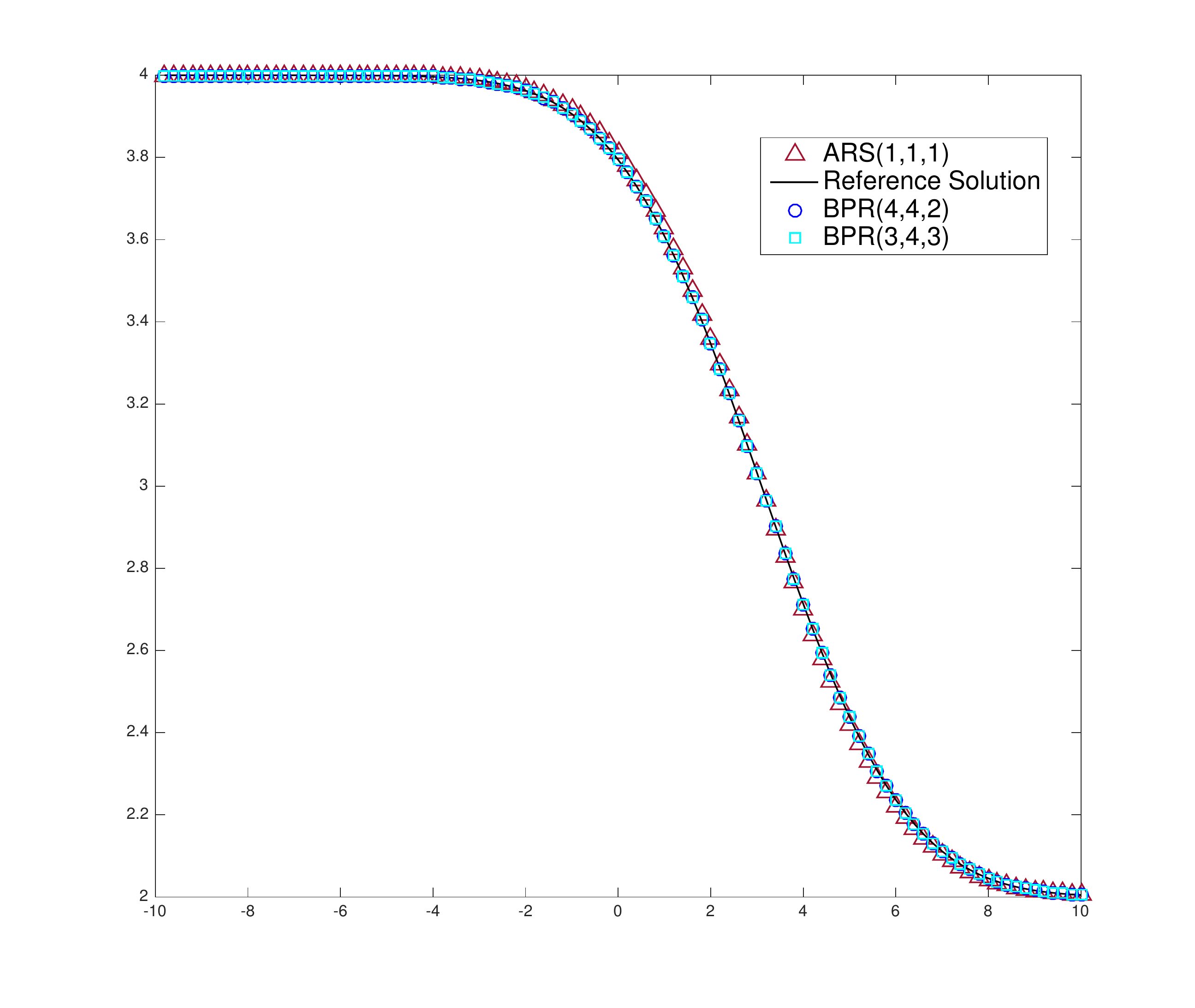}
\caption{Test 1. Solution of problem (\ref{E1}) with initial data (\ref{RWdata}), $\Delta x = 0.2$ and $\Delta t = 0.5 \Delta x$ for the density $u$. Left: the
rarefied regime $\varepsilon = 0.5$. Right: the parabolic regime $\varepsilon = 10^{-6}$.}
\label{LinearCase}
\end{center}
\end{figure}

\subsubsection*{Test 2. Multiscale limit in the nonlinear case}
Here we consider the nonlinear Ruijgrok-Wu model Eq.(\ref{E3}) for $\alpha \in [0,1]$ and interaction parameters  $c=2\epsi$ and $a=b=1$.

Adding and subtracting  the two equations in (\ref{E3})
one obtains the following macroscopic equations for $\rho$ and
$\j$
\be
\left\{  
\begin{array}{l} 
\displaystyle  
\partial_t \rho + \partial_x \j =0 \,, \\
\\
\displaystyle
\partial_t \j + \frac{1}{\epsi^{2\alpha}}\partial_x \rho
= {{1}\over{\epsi^{1+\alpha}}} \left\{- \j + \frac{1}{2}\left( \rho^2 - \epsi^2 \j^2\right) \right\}.
\\
\end{array}
\right. 
\label{E8}
\ee
For small values of $\epsi$ the model behaviour can be derived by the Chapman-Enskog expansion and is characterized by the viscous Burgers equation
\be
\begin{array}{l} 
\displaystyle \j = \frac{1}{2} \rho^2 -\varepsilon^{1-\alpha} \partial_x \rho
+ \varepsilon^{1+\alpha} \rho^2 \partial_x \rho+
{\cal O}(\epsi^2)\\
\displaystyle \partial_t \rho + \partial_x \left(\frac{\rho^2}{2}\right) = \epsi^{1+\alpha} \partial_{x} 
\left[\left(\frac1{\epsi^{2\alpha}}-\rho^2\right)\partial_x \rho \right] +{\cal O}(\epsi^2).
\end{array}
\ee 

We consider two different initial conditions. The first one is given by two local Maxwellian characterized by 
\be
\left\{
\begin{array}{ll}
\rho_L = 1.0,\quad \j_L=0,\quad & 
-10 < x < 0,
\\
\rho_R = 2.0,\quad \j_R=0,\quad &
0 < x < 10,
\\
\end{array} \right. \label{RWdata22}
\ee
with $\j = [(1+\rho^2\varepsilon^2)^{1/2} -1]/\varepsilon^2$. 

We show in Figure \ref{NonLinearCase} the numerical solution for $u$ in the case $\alpha=1$ using BPR(4,4,2) and BPR(3,4,3) schemes  in the rarefied ($\varepsilon = 0.4$) and parabolic ($\varepsilon = 10^{-6}$) regimes with $\Delta x = 0.2$ at the final time $T = 0.2$. The solution of both schemes is in very good agreement with the reference solution computed with $\Delta x=0.04$.

\begin{figure}[htb]
\begin{center}
\centering
\includegraphics[width=3.1in]{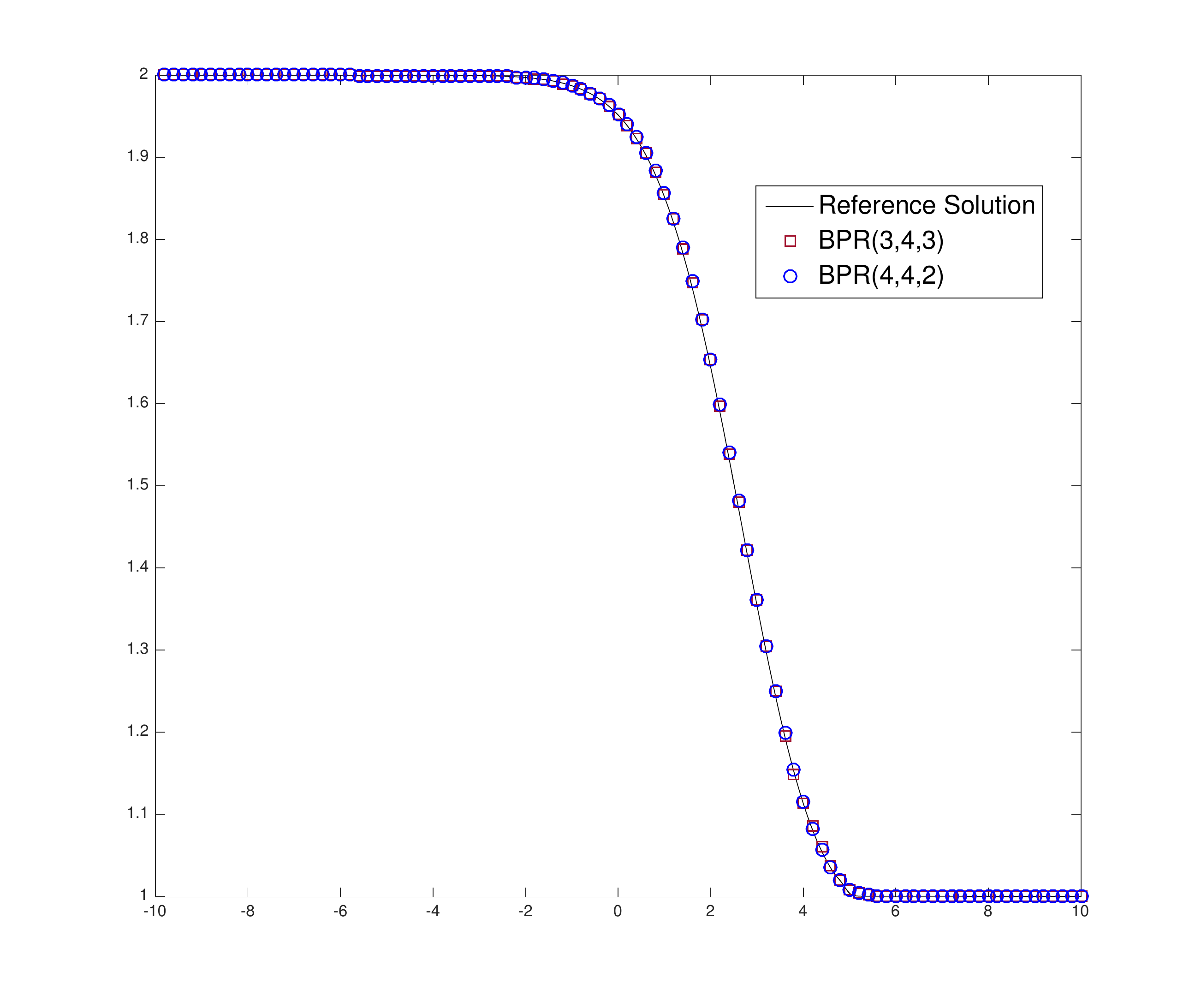},
\includegraphics[width=3.1in]{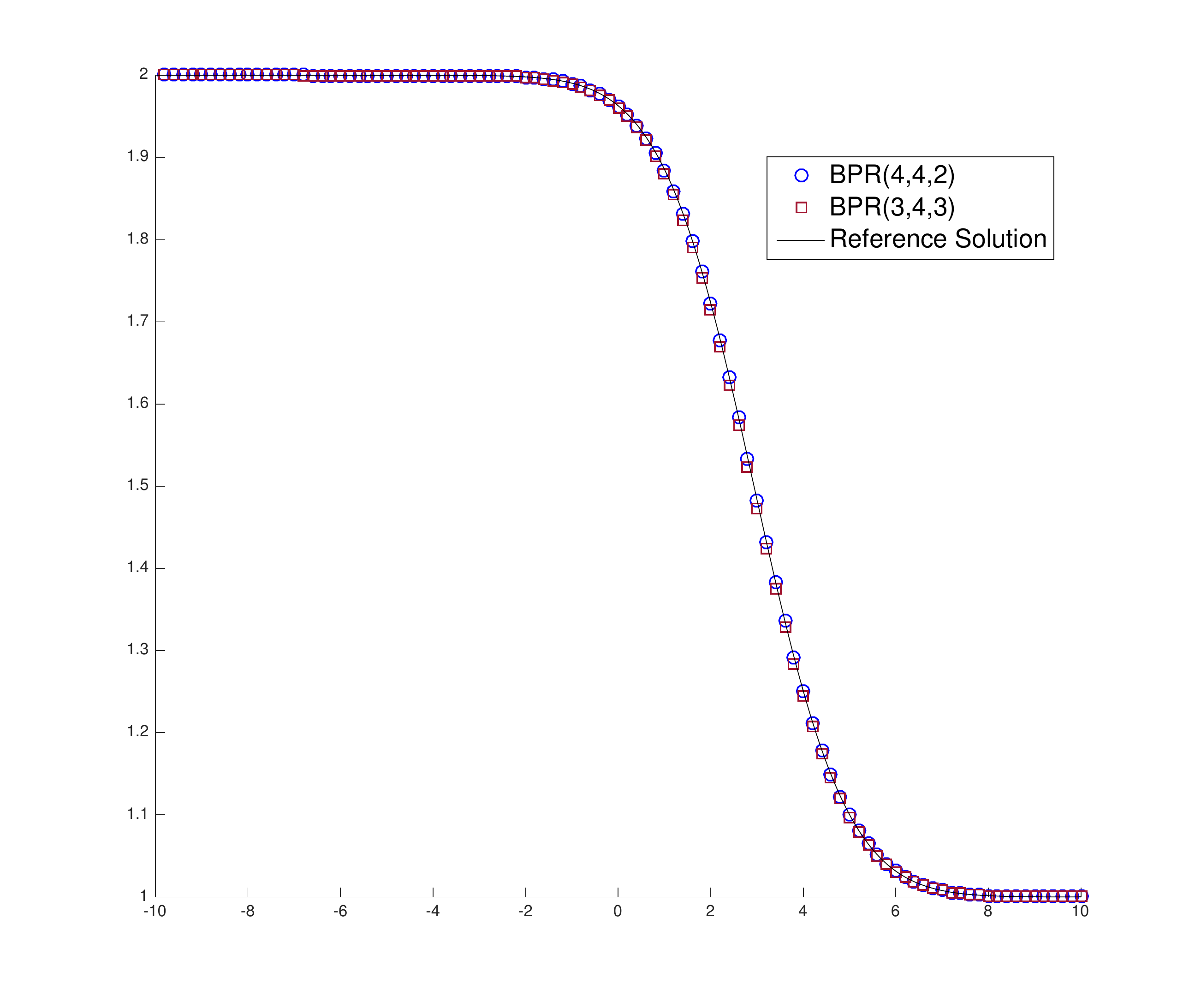}
\caption{Test 2. Numerical solutions of problem (\ref{E8}) for $\alpha=1$ with initial conditions (\ref{RWdata22}) at $T = 2.0$ and $\Delta t = 0.5 \Delta x$. Left: the rarefied regime
for $\rho$ with $\varepsilon = 0.4$. Right: the parabolic regime for $\rho$ with $\varepsilon  = 10^{-6}$. 
}
\label{NonLinearCase}
\end{center}
\end{figure}

\begin{figure}[htbp]
\begin{center}
\centering
\includegraphics[width=3.1in]{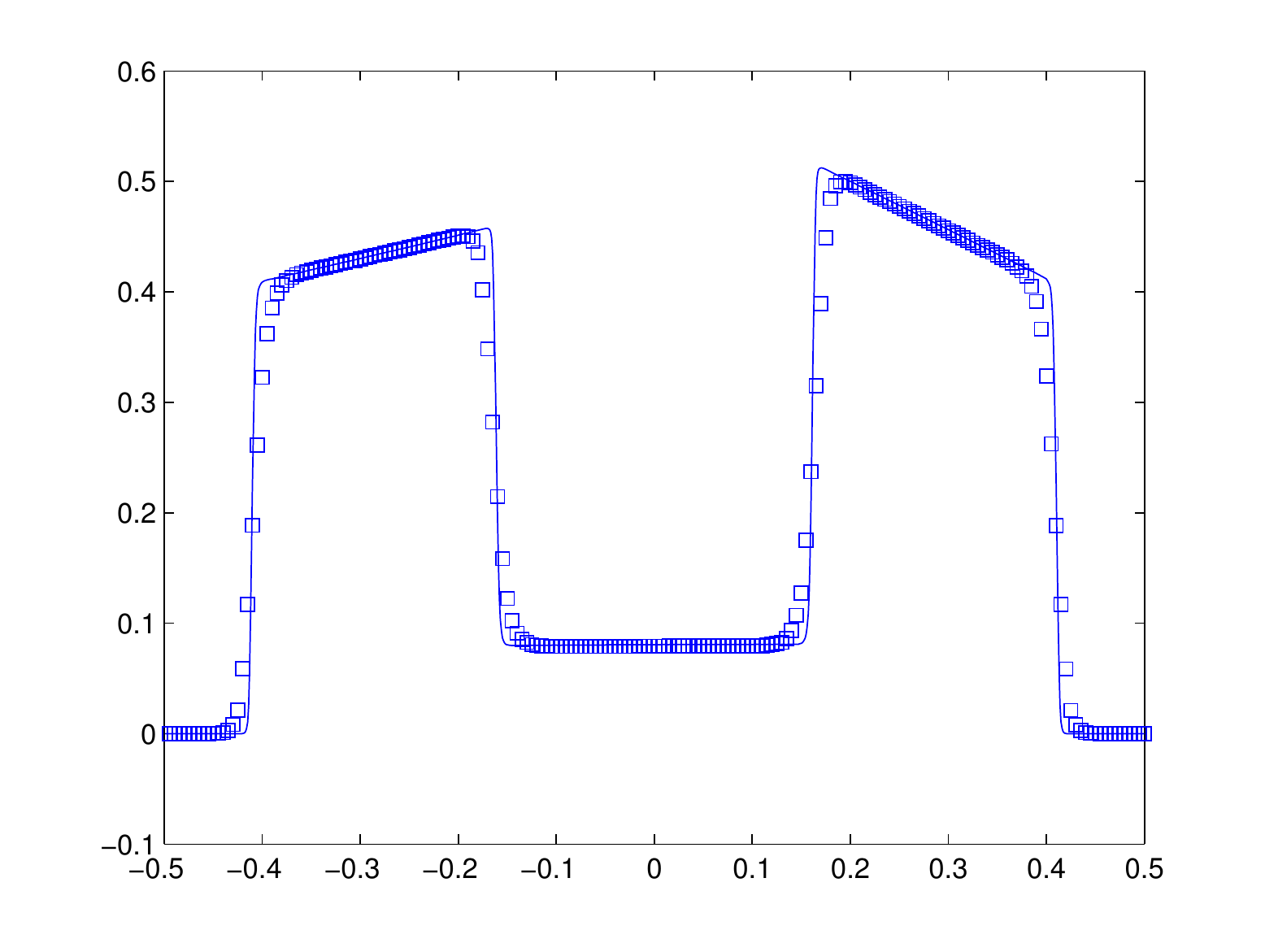},
\includegraphics[width=3.1in]{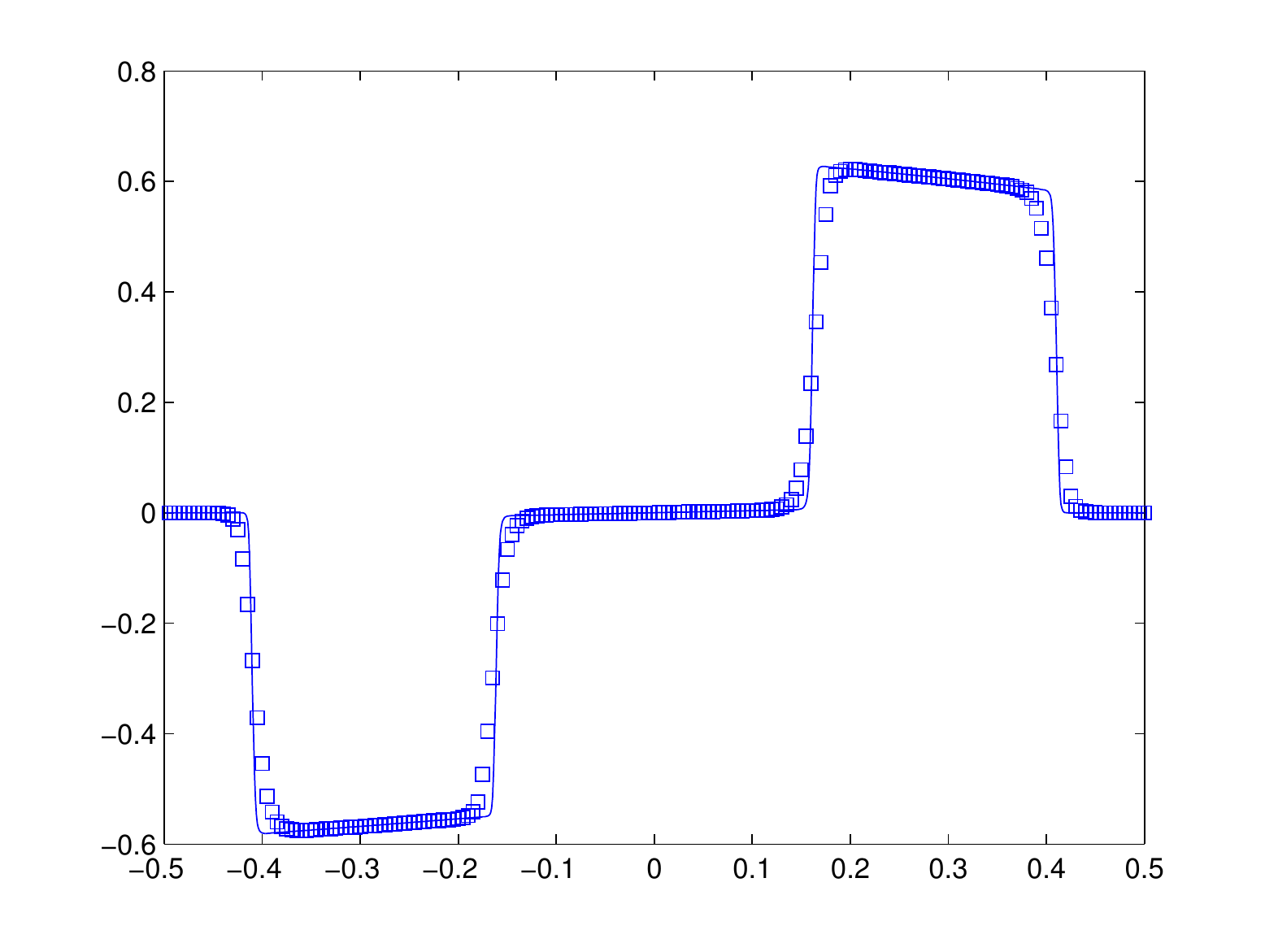}
\includegraphics[width=3.1in]{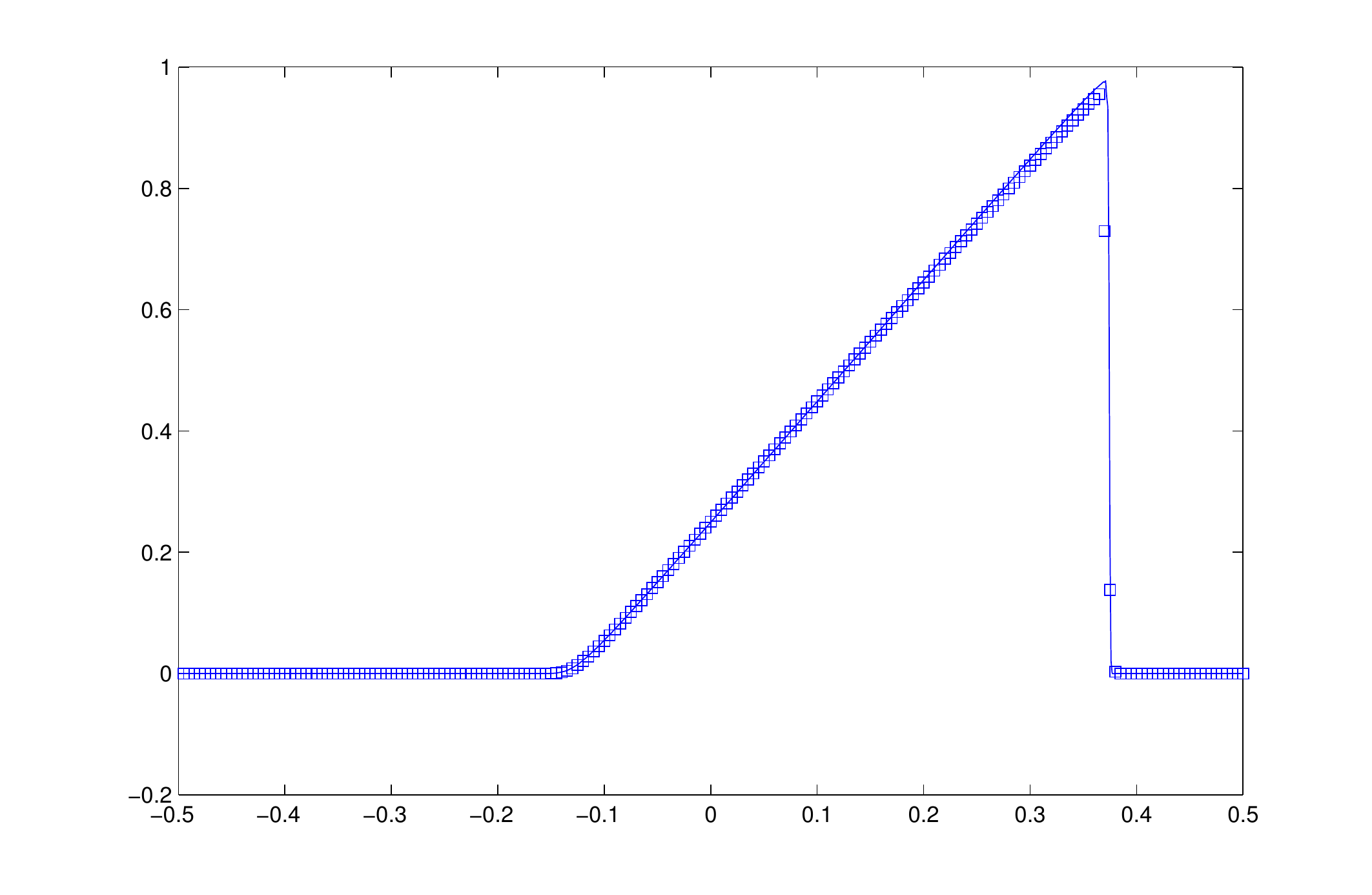}
\includegraphics[width=3.1in]{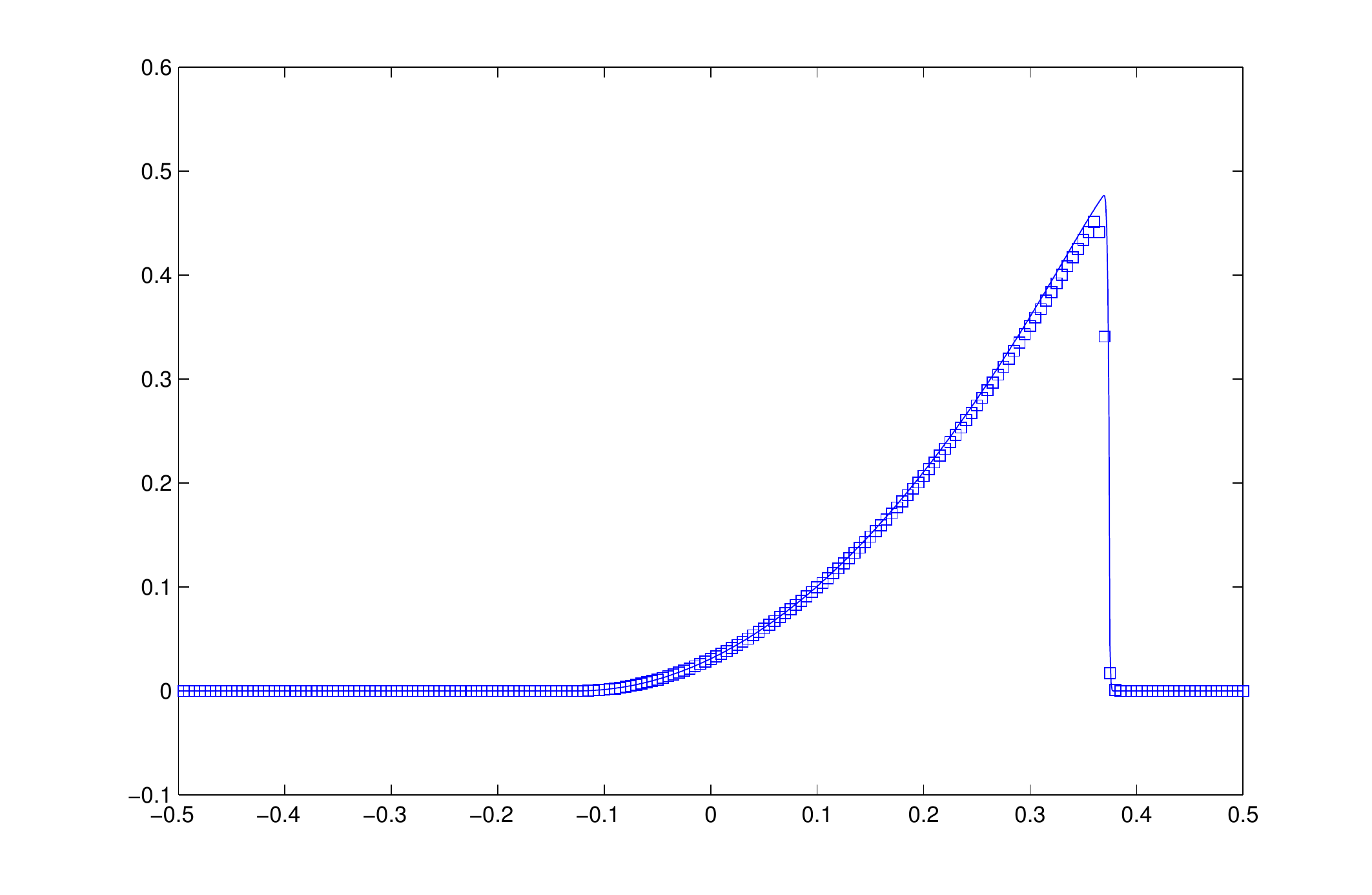}
\includegraphics[width=3.15in]{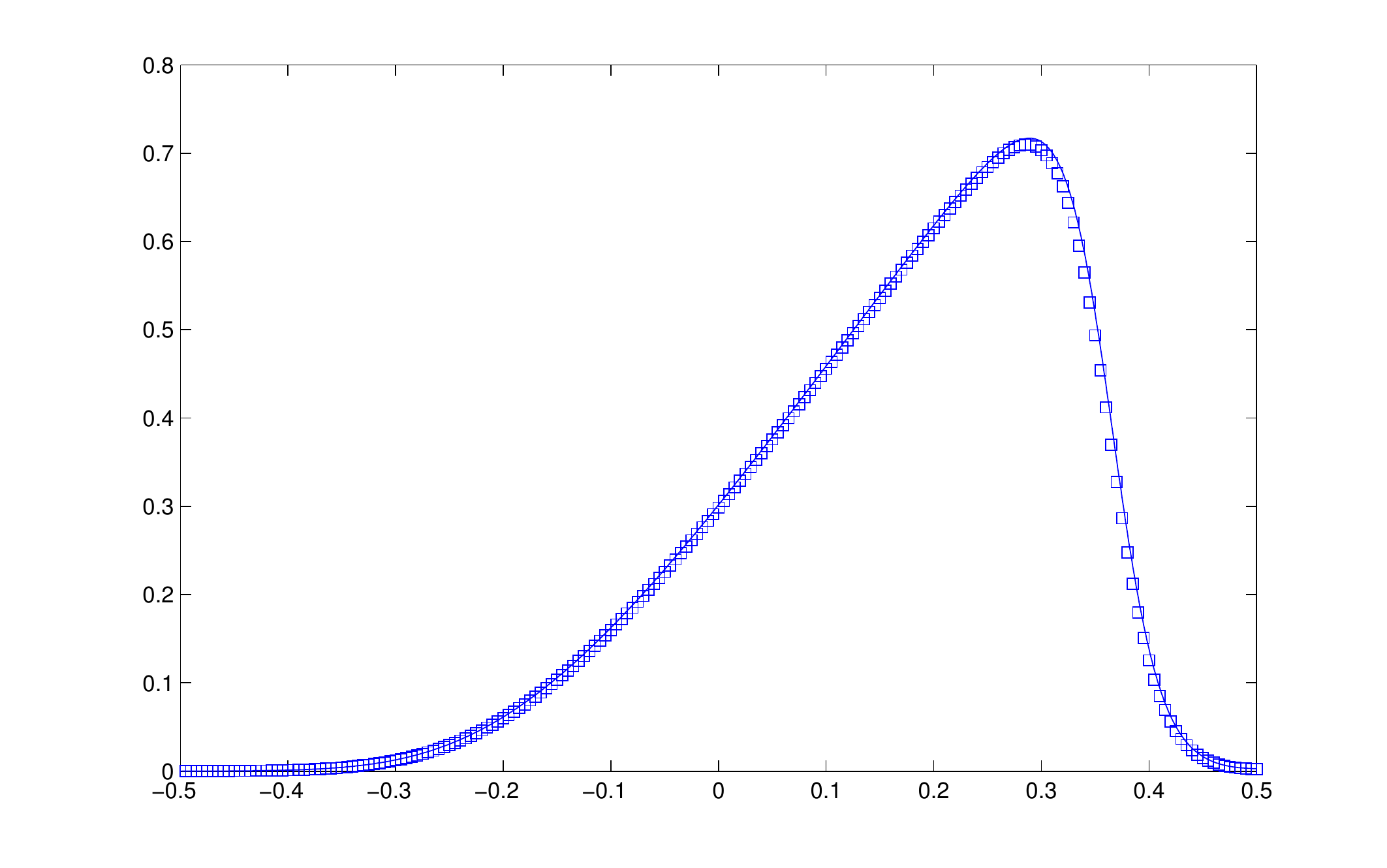}
\includegraphics[width=3.1in]{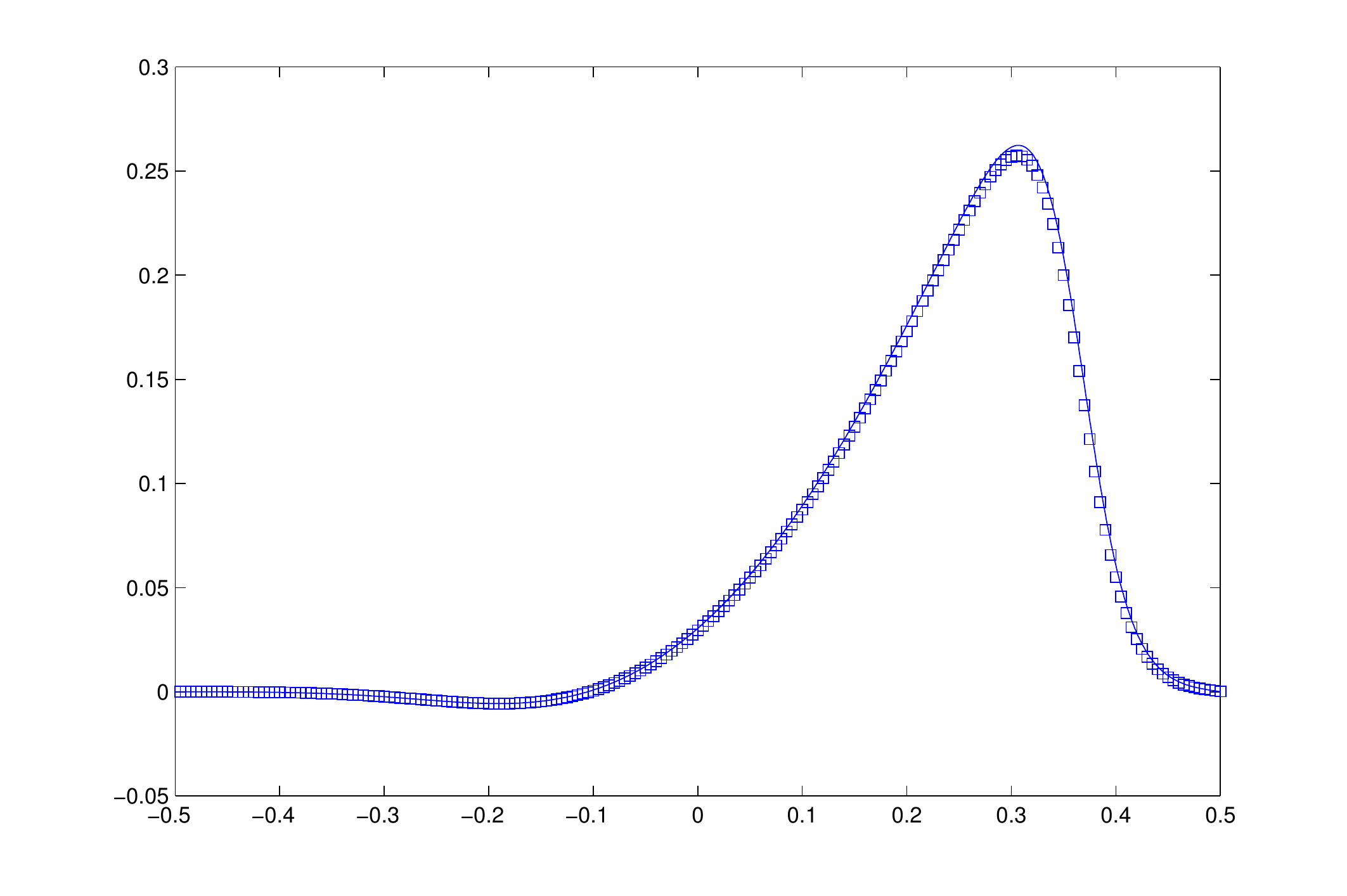}

\caption{Test 2. Numerical solutions of system (\ref{E8})  with initial conditions (\ref{riem}) for the mass density $\rho$ (left) and the momentum $\j$ (right) in the rarefied regime  (top panels) with $\varepsilon = 0.7$, $\alpha = 0$ and $\Delta t = 0.0025$, $\Delta x = 0.005$ at $T = 0.2$. In the parabolic regime $\varepsilon = 10^{-8}$, with $\alpha = 0.5$ (mid panels)  and $\alpha = 0.75$ (bottom panels) with $\Delta t = 0.004$ (i.e. $\Delta t = 0.8 \Delta x$) at time $T = 0.5$.}
\label{Case3}
\end{center}
\end{figure}

The last test case that we consider is the propagation of an initial square wave. The initial profile is specified as 
\be
\rho=1.0,\quad \j=0.0 \quad {\rm for}\quad |x| < 0.125,\qquad \rho = \j = 0 \quad {\rm for} \quad |x| > 0.125 \label{riem}
\ee
with reflecting boundary conditions, i.e. $v = 0$, $u_x = 0$ on the boundary.

We integrate the equations over $[-0.5, 0.5]$ with $200$ spatial cells. In Figure \ref{Case3} we plot the behavior of the system 
in the rarefied regime for $\varepsilon = 0.7$ with $\alpha = 0$  at time $T = 0.2$ and
in the parabolic regime ($\varepsilon = 10^{-8}$) with $\alpha = 0.5, 0.75$. We use BPR(3,4,3) scheme (\ref{ARS4}).
 The numerical solutions for the mass density $\rho$ and momentum $\j$ are computed in the rarefied regime and in the diffusive regime with  $\Delta t = 0.5\Delta x$ and are depicted with a {reference} solution obtained using fine grids with $\Delta x =  0.001$. We observe that in the rarefied regime no oscillation
appears and the scheme describes well the behavior of the system. In the parabolic regime the choice of the parameter $\alpha$ gives  the corresponding Reynolds numbers of the problem, i.e., $Re = 10000$ and $Re = 100$. The output of the solution in the parabolic regime is given at $T = 0.5$ with $\Delta t = 0.004$ (i.e. $\Delta t = 0.8 \Delta x$) and a right moving shock is formed. Although $\varepsilon$ is underresolved the scheme proposed captures well the correct behavior of the equilibrium equation independently of $\alpha$.

\subsubsection*{Test 3. Multiscale space varying limit in the nonlinear case}
Finally, we consider the multi scale relaxation system of the previous section, in the case where the parameter $\alpha=\alpha(x)\in [0, 1]$ depends on the space variable. Therefore, the limit behavior of the system may depend on the particular region of the computational domain. 

Now we apply our schemes to system (\ref{E8}) considering two different cases of the varying $\alpha$ number. The domain is chosen to be $x \in [-0.5, 0.5]$. In the first case $\alpha$  increases smoothly from a small value $\alpha_0$  to $\mathcal{O}(1)$, by the formula,
\be\label{alpha1}
\alpha(x) = \alpha_0 + 0.5(1 + \tanh(20(x+ 0.1)))
\ee
with $\alpha_0 = 10^{-6}$.
In the second case we consider  the function $\alpha$ which contains a discontinuity
\be\label{alpha2}
\left\{
\begin{array}{l}
\alpha_L = 0.0,\quad 
-0.5 < x < 0,
\\
\alpha_R = 1.0,\quad 
0 < x < 0.5.
\\
\end{array} \right. 
\ee
The two different scenarios for $\alpha$ are depicted in the Figure \ref{alpha}. 

\begin{figure}[htbp]
\begin{center}
\centering
\includegraphics[width=3.1in]{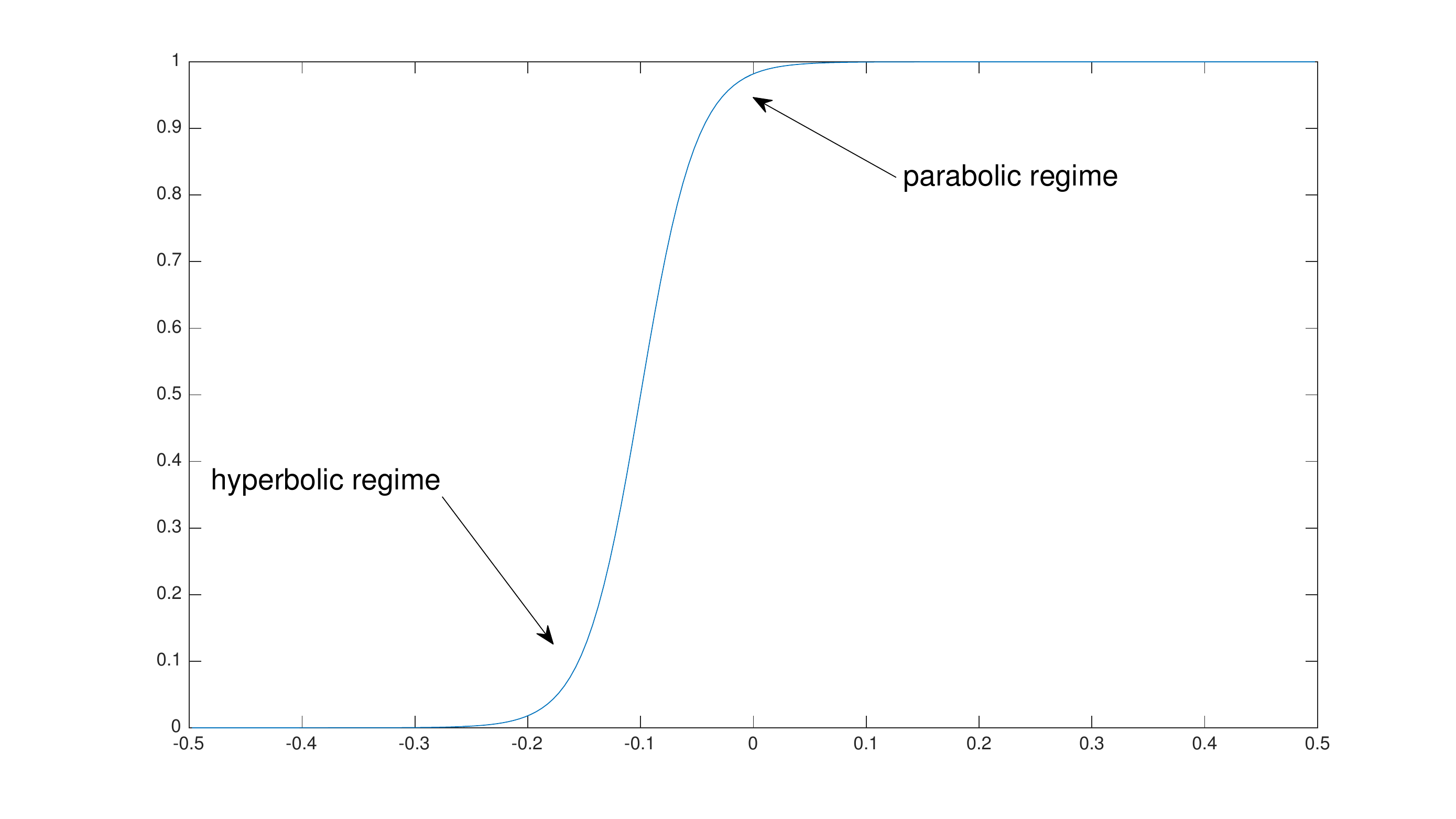},
\includegraphics[width=3.1in]{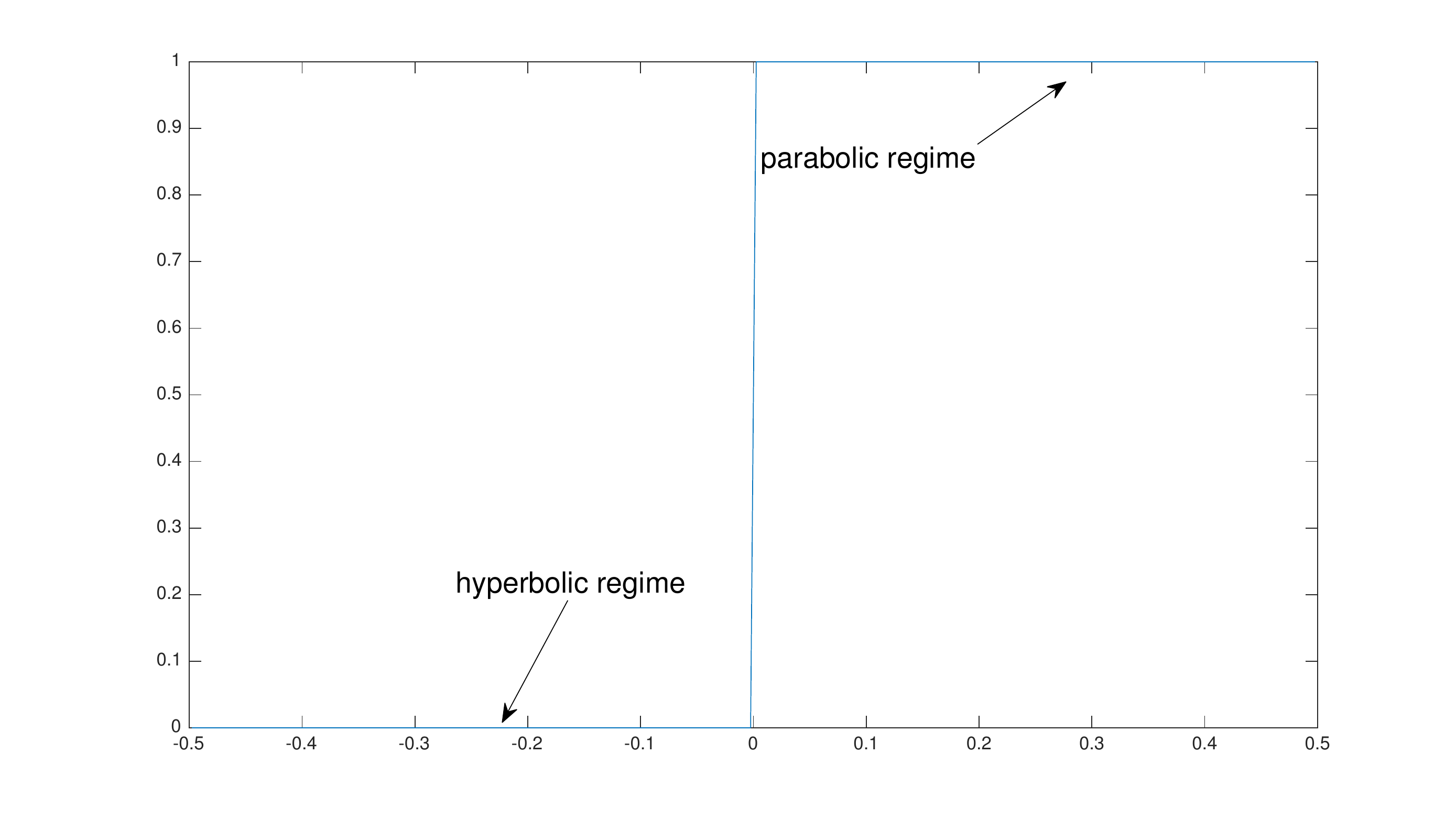}
\caption{Test 3. Space varying $\alpha(x)$. Left: the smooth profile (\ref{alpha1}). Right: the discontinouos profile (\ref{alpha2})}
\label{alpha}
\end{center}
\end{figure}

The initial profile is specified by (\ref{riem}) with reflecting boundary conditions and we fix $\varepsilon = 10^{-8}$ and $\Delta t = 0.5 \Delta x$. Here we use BPR(3,3,5) scheme. The reference solutions are computed using a fine grid with $\Delta x = 0.001$. The results are depicted up to time $T = 0.05$ (left) and $T = 0.18$ (right).  In Figure \ref{alphaFinal} we observe that the scheme is able to capture the correct behavior of the reference solution even in this test case.

\begin{figure}[htbp]
\begin{center}
\centering
\includegraphics[width=3.1in]{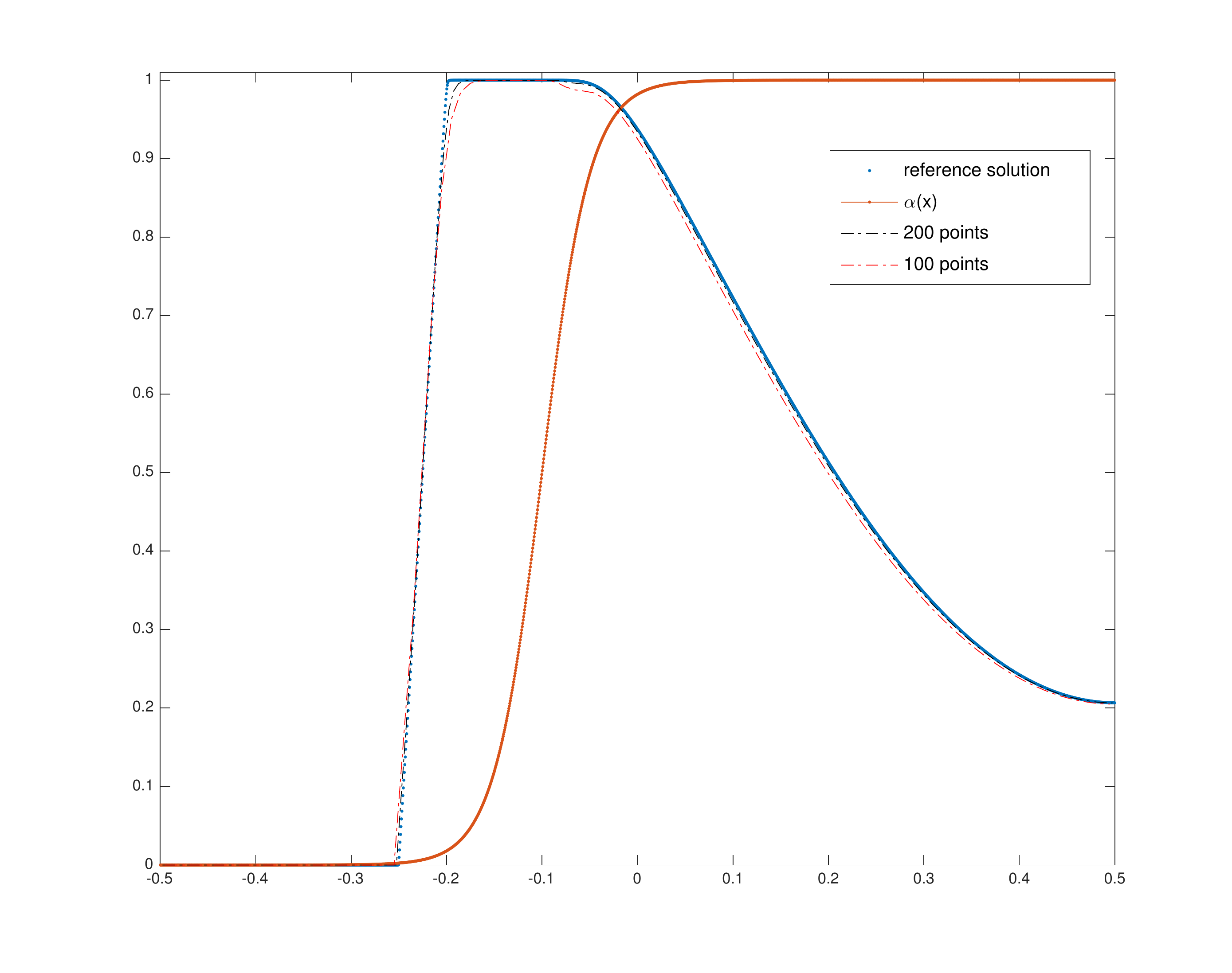}
\includegraphics[width=3.1in]{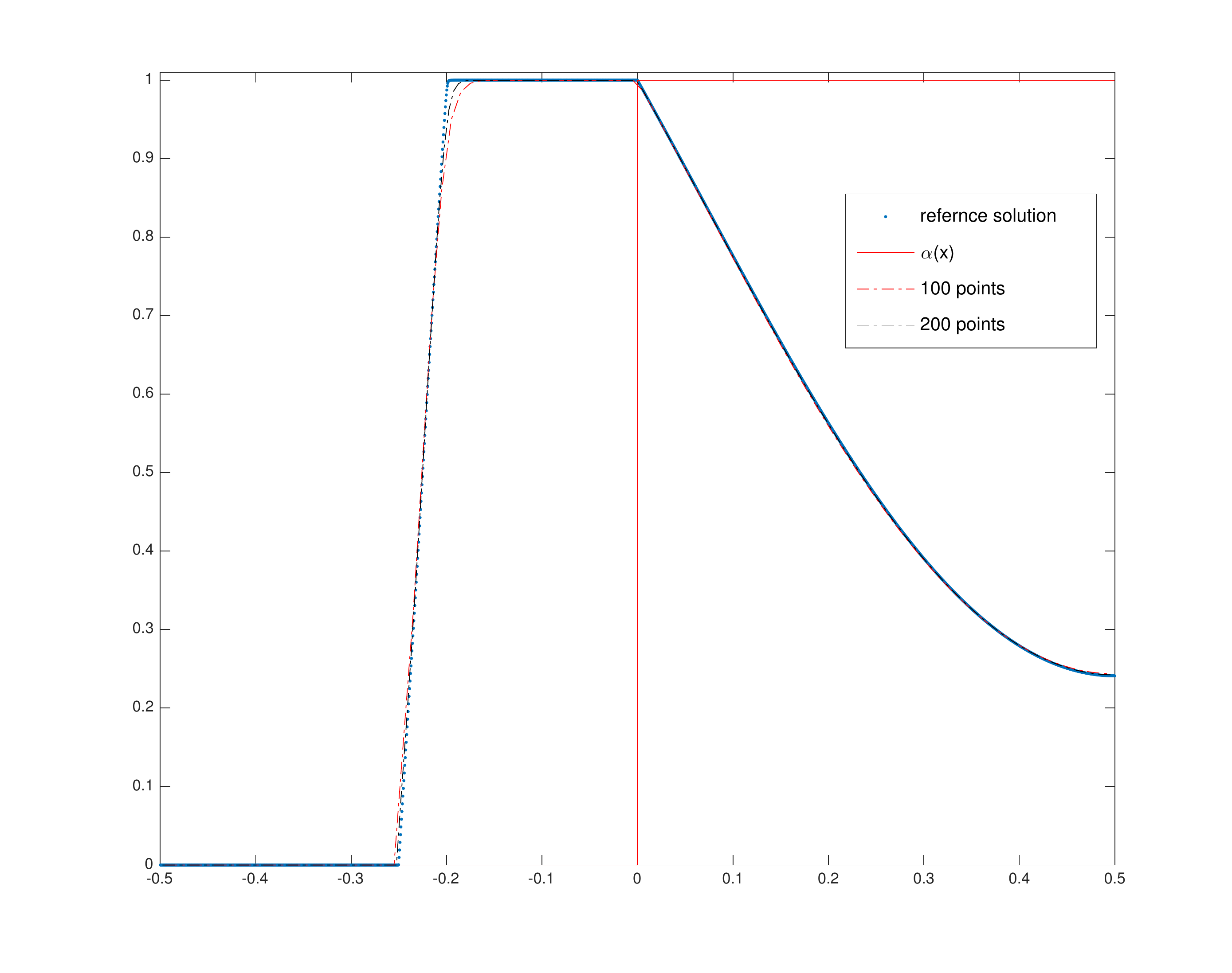}
\caption{Test 3. For $\varepsilon = 10^{-8}$, the left and right pictures show the behaviour of the component $u$ in the two different cases of $\alpha$.}
\label{alphaFinal}
\end{center}
\end{figure}

\section{Conclusions} 
In this work we have developed a new IMEX-RK approach that is capable correctly to capture the asymptotic behavior of hyperbolic balance laws with relaxation under different kinds of scaling. Previous IMEX-RK schemes were designed to deal specifically with one kind of scaling, either hyperbolic \cite{BR2009, CJR, Jin, PR} or parabolic \cite{BPR, Klar, JPT, NP}, whereas the present schemes are robust enough to be able to deal with both scalings. Related approaches were presented in \cite{BR,JP, Lem2, NP2}.

From a physical viewpoint, these scaling limits corresponds to the classical fluid-dynamic scalings in the kinetic theory of rarefied gases that lead from the Boltzmann equation to its compressible and incompressible limits. Several numerical test for a simple kinetic model which possesses different asymptotic limits have confirmed the validity of the present approach. In the near future, we hope to extend this class of IMEX-RK schemes to the full Boltzmann equation by adopting the penalization techniques developed in \cite{DP, JF}.  

\subsection*{Acknowledgments} The authors are grateful to Prof. M. Lemou and Prof. G. Dimarco for helpful observations. The research that led to the present paper was partially supported by the research grant {\it Numerical methods for uncertainty quantification in hyperbolic and kinetic equations} 
of the group GNCS of INdAM, by 
ITN-ETN Marie-Curie Horizon 2020 program ModCompShock, {\em Modeling and computation of shocks and interfaces}, Project ID: 642768 and by the INDAM-GNCS 2017 research grant {\em Numerical methods for hyperbolic and kinetic equation and applications}.

\appendix
\section{Appendix: IMEX-RK schemes}
{
In order to achieve higher then first order accuracy in time for (GSA) IMEX-RK schemes, we need to observe first that there are no second order GSA IMEX-RK methods of type I with $s= 3$, i.e. with three stages \cite{BPR}. Thus, to construct a second order GSA IMEX Runge-Kutta schemes, we consider the following proposition \cite{BPR} 
\begin{proposition} 
The only type of second order GSA IMEX-RK scheme with three levels $s = 3$ 
that satisfies the classical second order conditions (\ref{OrdCond}) is the type II  where $c_i = \tilde{c}_i$ for all $i = 2,...,s-1$, with $c_1 = 0$ and $c_s = 1$.
\end{proposition}
}

{
As usual, order conditions are obtained by matching the Taylor expansion of the
exact solution and the numerical one, up to terms of the prescribed order. In the
relaxed case, i.e. $\varepsilon \to 0$, system  (\ref{I61}) for $\alpha=1$ reduces to
\begin{equation}\label{App1} v = -p(u)_x + f(u), \quad u_t + f(u)_x = p(u)_{xx}.\end{equation}
The unified IMEX-RK approach described in Section \ref{SecPar} provides a scheme that converges to an explicit-implicit scheme for the limit convection-diffusion equation in (\ref{App1}). Performing the same analysis as proposed in \cite{BR} for system  (\ref{I61}) we obtain for the $u$-component the classical order conditions, while some additional order conditions for the $v$-component are required in order to have consistency and maintain the classical order  in the limit case.  We recall that the GSA  assumption of the method guarantees that IMEX-RK scheme relaxes at the same IMEX-RK one  when $\varepsilon \to 0$, but in order to maintain the order of accuracy of the scheme in the limit we must impose some additional order conditions.}

{Below we list these new additional order conditions that the $v$-component must satisfy up to second order
\begin{eqnarray}\label{AddOrdCondNew}
\begin{array}{|l|llll|}
\hline&&&&\\[-.2cm]
\mathrm{consistency} & b^TA^{-2}\tA e = 1,&&&\\
\mathrm{first \ order} & b^TA^{-2}\tA c = 1,& b^TA^{-2}\tA \tilde{c} = 1,&&\\
\mathrm{second\ order} & b^TA^{-2}\tA c^2 = 1, & b^TA^{-2}\tA \tilde{c}^2 = 1, & b^TA^{-2}\tA \tilde{c}c = 1, & b^TA^{-2}\tA c\tilde{c} = 1,\\
& b^TA^{-2}\tA A c = 1/2,& b^TA^{-2}\tA A \tilde{c} = 1/2, & &\\
& b^TA^{-2}\tA \tA c = 1/2, & b^TA^{-2}\tA \tA \tilde{c} = 1/2.&&\\[+.2cm]
\hline
\end{array}
\end{eqnarray}
Note that in order to reduce the number of the additional order conditions we require that $\tilde{c} = c$, this assumption simplifies a lot the number of coupling conditions 
\begin{eqnarray}\label{AddOrdCond2}
\begin{array}{|l|lll|}
\hline&&&\\[-.2cm]
\mathrm{consistency}& b^TA^{-2}\tA e = 1,&&\\
\mathrm{first \ order}& b^TA^{-2}\tA c = 1,&&\\
\mathrm{second\ order}& b^TA^{-2}\tA c^2 = 1, & b^TA^{-2}\tA A c = 1/2, & b^TA^{-2}\tA \tA c = 1/2.\\[+.2cm]
\hline
\end{array}
\end{eqnarray}
}

}  

Finally, we present the different IMEX-RK schemes, up to order three, used in our numerical tests. 
Below, these schemes are represented as usual by the double Butcher tableau. On the left we have the explicit part and on the right the implicit part of the IMEX schemes. All the schemes satisfy the GSA property, but only the new scheme BPR$(4,4,2)$ satisfies the additional order conditions (\ref{AddOrdCond2}).

\begin{itemize}
\item First order ARS$(1,1,1)$ scheme  
\begin{eqnarray}\label{ars1}
\begin{array}{c|cc}
              0  & 0 &0 \\
              1 & 1 & 0  \\
              \hline 
                 & 1 & 0 
\end{array} \qquad
\begin{array}{c|cc}
              0  & 0 & 0 \\
              1  & 0 & 1 \\
              \hline
                 & 0 & 1  
\end{array}
\end{eqnarray}

\item Second order ARS$(2,2,2)$ scheme \cite{ARS} 
\begin{eqnarray}\label{ars2}
\begin{array}{c|ccc}
              0  & 0 & 0 &0 \\
              \gamma & \gamma & 0 &0  \\
              1 & \delta & 1-\delta & 0\\
              \hline 
                 & \delta & 1-\delta & 0 
\end{array} \qquad
\begin{array}{c|ccc}
              0  & 0 & 0 & 0  \\
              \gamma  & 0 & \gamma & 0 \\
              1 & 0 & 1-\gamma & \gamma \\
              \hline
                 & 0 & 1-\gamma & \gamma 
\end{array}.
\end{eqnarray}
where $\gamma = 1-1/\sqrt{2}$ and $\delta= 1- 1/2\gamma$, 
 \item Second order CK$(2,2,2)$ scheme \cite{CK}
 \begin{eqnarray}\label{PtypeCK}
\begin{array}{c|ccc}
              0 & 0 & 0 & 0\\
              2/3 &  2/3 & 0&0\\
              1 & 1/4 & 3/4 & 0\\
              \hline
              & 1/4  & 3/4 & 0
\end{array} \qquad
\begin{array}{c|ccc}
              0 & 0  & 0 & 0\\
               2/3 & -1/3+\sqrt{2}/2 & 1-\sqrt{2}/2 & 0\\
              1 & 3/4-\sqrt{2}/4& -3/4+3\sqrt{2}/4 & 1-\sqrt{2}/2 \\ 
              \hline 
              &3/4-\sqrt{2}/4 & -3/4+3\sqrt{2}/4 & 1-\sqrt{2}/2
\end{array}.
\end{eqnarray}

\item Second order BPR$(4,4,2)$ scheme

\begin{eqnarray} \label{New}
\begin{array}{c|ccccc}
0& 0& 0& 0& 0&0\\
1/4& 1/4& 0& 0&0&0\\
1/4&13/4&-3& 0& 0&0\\
3/4& 1/4&0& 1/2& 0& 0\\
1& 0&1/3& 1/6& 1/2& 0\\
\hline
 & 0&1/3 & 1/6& 1/2 & 0\\
\end{array} \ \ \ \ \ 
\begin{array}{c|ccccc}
0& 0& 0& 0& 0 & 0\\
1/4&0& 1/4& 0& 0 &0\\
1/4&0&0& 1/4& 0 & 0\\
3/4& 0& 1/24& 11/24& 1/4& 0\\
1& 0 &11/24 & 1/6 & 1/8& 1/4\\
\hline
& 0&11/24& 1/6& 1/8& 1/4\\
\end{array}.
\end{eqnarray}

\item Third order BPR$(3,4,3)$ scheme \cite{BPR}
\begin{eqnarray} \label{ARS4}
\begin{array}{c|ccccc}
0& 0& 0& 0& 0&0\\
1& 1& 0& 0&0&0\\
2/3&4/9&2/9& 0& 0&0\\
1& 1/4&0& 3/4& 0& 0\\
1& 1/4&0& 3/4& 0& 0\\
\hline
 & 1/4&0 & 3/4& 0 & 0\\
\end{array} \ \ \ \ \ 
\begin{array}{c|ccccc}
0& 0& 0& 0& 0 & 0\\
1&1/2& 1/2& 0& 0 &0\\
2/3&5/18&-1/9& 1/2& 0 & 0\\
1& 1/2& 0& 0& 1/2& 0\\
1& 1/4 & 0 & 3/4& -1/2& 1/2\\
\hline
& 1/4& 0& 3/4& -1/2& 1/2\\
\end{array}.
\end{eqnarray}

\end{itemize}




\end{document}